\documentclass[reqno,11pt]{amsart}

\usepackage{amssymb,amsmath,amsthm,amsfonts,mathtools,bbm,mathrsfs,enumitem,graphics,float,multirow,enumitem,cases}
\usepackage[normalem]{ulem}
\usepackage{hyperref}
\usepackage{comment}
\usepackage{hyperref}
\hypersetup{
    colorlinks,
    linkcolor={blue!80!black},
    citecolor={red!80!black},
    urlcolor={red!80!black}
}
\usepackage{xargs}                      
\usepackage{xcolor}  
\usepackage[colorinlistoftodos,prependcaption,textsize=tiny]{todonotes}
\newcommandx{\unsure}[2][1=]{\todo[linecolor=red,backgroundcolor=red!25,bordercolor=red,#1]{#2}}
\newcommandx{\change}[2][1=]{\todo[linecolor=blue,backgroundcolor=orange!25,bordercolor=blue,#1]{#2}}
\newcommandx{\info}[2][1=]{\todo[linecolor=OliveGreen,backgroundcolor=OliveGreen!25,bordercolor=OliveGreen,#1]{#2}}

\usepackage[top=1.6in,bottom=1.5in,left=1.3in,right=1.3in]{geometry}
\allowdisplaybreaks
\usepackage{cellspace}
\newcommand\redout{\bgroup\markoverwith{\textcolor{red}{\rule[0.5ex]{2pt}{0.8pt}}}\ULon}

\newtheorem{theorem}{Theorem}
\newtheorem{lemma}[theorem]{Lemma}
\newtheorem{claim}[theorem]{Claim}
\newtheorem{fact}[theorem]{Fact}

\newtheorem{definition}{Definition}

\usepackage{tikz}
\usetikzlibrary{calc}

\tikzset{font= }

\newcommand\nc\newcommand
\nc\bfa{{\boldsymbol a}}\nc\bfA{{\boldsymbol A}}\nc\cA{{\mathscr A}}
\nc\bfb{{\boldsymbol b}}\nc\bfB{{\boldsymbol B}}\nc\cB{{\mathscr B}}
\nc\bfc{{\boldsymbol c}}\nc\bfC{{\boldsymbol C}}\nc\cC{{\mathscr C}}
\nc\bfd{{\boldsymbol d}}\nc\bfD{{\boldsymbol D}}\nc\cD{{\mathscr D}}
\nc\bfe{{\boldsymbol e}}\nc\bfE{{\boldsymbol E}}\nc\cE{{\mathscr E}}
\nc\bff{{\boldsymbol f}}\nc\bfF{{\boldsymbol F}}\nc\cF{{\mathscr F}}
\nc\bfg{{\boldsymbol g}}\nc\bfG{{\boldsymbol G}}\nc\cG{{\mathscr G}}
\nc\bfh{{\boldsymbol h}}\nc\bfH{{\boldsymbol H}}\nc\cH{{\mathscr H}}\nc\fH{{\mathfrak H}}
\nc\bfi{{\boldsymbol i}}\nc\bfI{{\boldsymbol I}}\nc\cI{{\mathcal I}}
\nc\bfj{{\boldsymbol j}}\nc\bfJ{{\boldsymbol J}}\nc\cJ{{\mathscr J}}
\nc\bfk{{\boldsymbol k}}\nc\bfK{{\boldsymbol K}}\nc\cK{{\mathscr K}}
\nc\bfl{{\boldsymbol l}}\nc\bfL{{\boldsymbol L}}\nc\cL{{\mathscr L}}
\nc\bfm{{\boldsymbol m}}\nc\bfM{{\boldsymbol M}}\nc\cM{{\mathscr M}}
\nc\bfn{{\boldsymbol n}}\nc\bfN{{\boldsymbol N}}\nc\sN{{\mathscr N}}
\nc\bfo{{\boldsymbol o}}\nc\bfO{{\boldsymbol O}}\nc\cO{{\mathscr O}}
\nc\bfp{{\boldsymbol p}}\nc\bfP{{\boldsymbol P}}\nc\cP{{\mathscr P}}
\nc\bfq{{\boldsymbol q}}\nc\bfQ{{\boldsymbol Q}}\nc\cQ{{\mathscr Q}}
\nc\bfr{{\boldsymbol r}}\nc\bfR{{\boldsymbol R}}\nc\cR{{\mathscr R}}
\nc\bfs{{\boldsymbol s}}\nc\bfS{{\boldsymbol S}}\nc\cS{{\mathscr S}}
\nc\bft{{\boldsymbol t}}\nc\bfT{{\boldsymbol T}}\nc\cT{{\mathscr T}}
\nc\bfu{{\boldsymbol u}}\nc\bfU{{\boldsymbol U}}\nc\cU{{\mathscr U}}
\nc\bfv{{\boldsymbol v}}\nc\bfV{{\boldsymbol V}}\nc\cV{{\mathscr V}}
\nc\bfw{{\boldsymbol w}}\nc\bfW{{\boldsymbol W}}\nc\cW{{\mathscr W}}
\nc\bfx{{\boldsymbol x}}\nc\bfX{{\boldsymbol X}}\nc\cX{{\mathscr X}}
\nc\bfy{{\boldsymbol y}}\nc\bfY{{\boldsymbol Y}}\nc\cY{{\mathscr Y}}
\nc\bfz{{\boldsymbol z}}\nc\bfZ{{\boldsymbol Z}}\nc\cZ{{\mathscr Z}}
\nc\pp{\mathbb{P}}
\nc\ee{\mathbb{E} }

\renewcommand{\le}{\leqslant}
\renewcommand{\leq}{\leqslant}
\renewcommand{\ge}{\geqslant}
\renewcommand{\geq}{\geqslant}

\DeclareMathOperator{\Var}{Var}
\nc{\Cay}{{\sf Cay}}

\nc{\ff}{{\mathbb F}}

\title[]{Bootstrap percolation on a generalized Hamming cube \MakeUppercase{\romannumeral 2}}
\author[]{Fengxing Zhu}\thanks{Institute for Systems Research and Department of ECE, University of Maryland, College Park, MD 20742, USA, fengxing@terpmail.umd.edu. Supported in part by NSF grant CCF 2330909.}
\date{}

\begin{document}
\begin{abstract}
    In this paper we obtain the critical probability  $p_c(Q_{k,n},r)$ for bootstrap percolation with the infection threshold $r=\frac{N}{2}$ on the generalized $n$-dimensional hypercube $Q_{k,n}$ with vertex set $V(Q_{k,n})=\{0,1\}^n$ and edges connecting the pairs at Hamming distance $1,2,\dots,k$, where $k\ge 2$ and $N=\sum_{i=1}^k\binom{n}{i}$. More precisely, we obtain the exact first-order and second-order terms of $p_c(Q_{k,n},\frac{N}{2})$ by extending the main theorem of Balogh, Bollob{\'a}s, and Morris (2009). 
\end{abstract}
\maketitle

\section{Introduction}
The process of $r$-neighbor bootstrap percolation on an undirected graph $G(V,E)$, with an integer $r\geq 1$, was introduced by Chalupa, Leith, and Reich \cite{Chalupa_1979}. In this process, each vertex is either infected or healthy, and once a vertex becomes infected, it remains infected forever. Initially, a set of vertices $A_0$ is infected, and let $A_i$ denote the set of infected vertices up to step $i$. The bootstrap percolation process evolves in discrete steps as follows: for $i>0$,
$$A_i = A_{i-1} \cup \{v \in V:|\mathcal{N}(v) \cap A_{i-1}| \geq r \},$$
where $\mathcal{N}(v)$ is the neighborhood of the vertex $v$ in $G$. In simple terms, a vertex that is not initially infected becomes infected at step $i$ if it has at least $r$ infected neighbors at step $i-1$. Here, we define $r$ as the infection threshold for all $v \in V$ and a contagious set as a set of initially infected vertices that leads to the complete infection of the entire graph. Percolation occurs if all vertices have been infected by the end of the process. In random bootstrap percolation, each vertex is independently and randomly infected at the beginning with a probability $p$. The central question in random bootstrap percolation is to determine the critical probability,
$$p_c(G,r):=\sup\{p \in (0,1): \mathbb{P}_p(A \; \text{percolates on} \: G) \leq \frac{1}{2} \},$$
where $A$ represents the set of initially infected vertices.

For bootstrap percolation in the random setting, Aizenman and Lebowitz \cite{Aizenman_1988} demonstrated that when the process occurs within a finite box $\{1, 2, \ldots, n\}^d$, the critical probability scales as ${\Theta}((\log(n))^{1-d})$ with $d \geq r=2$. In their work \cite{BALOGH2003305}, Balogh and Bollob{\'a}s established sharp transition for bootstrap percolation on a $d$-dimensional box, building upon the results from \cite{Friedgut}. Holroyd \cite{Holroyd2002SharpMT} achieved a significant breakthrough by identifying the precise threshold for $d=r=2$. Subsequent advancements were made by Gravner and Holroyd \cite{holroyd_2008} as well as by Hartarsky and Morris \cite{Robert}, who determined the scaling of the second-order term,
    $$p_c([n]^2,r=2)=\frac{\pi}{18 \log n }-\frac{\Theta(1)}{(\log n)^{\frac{3}{2}}}.$$
   For the case where $d=r=3$, Cerf and Cirillo \cite{Cirillo} established that the critical probability scales as ${\Theta} \left(\frac{1}{\log \log(n)} \right)$. For $d \geq r \geq 3$, Cerf and Manzo \cite{CERF200269} demonstrated that the critical probability behaves as  ${\Theta}\left((\log^{*(r-1)}(n))^{r-d-1} \right)$, where $\log^{*(r-1)}$ represents the logarithm iterated $r-1$ times. The sharp threshold for the $r$-neighbor process on $[n]^d$ was ultimately determined by Balogh, Bollob{\'a}s, Duminil-Copin, and Morris \cite{Balogh}, extending the results by Balogh, Bollob{\'a}s, and Morris \cite{Balogh2}. For every $d \geq r \geq 2$, they identified an explicit constant $\lambda(d,r) > 0$ such that the critical probability can be expressed as follows: 
   $$
   p_c([n]^d,r)=\left[ \frac{\lambda(d,r)+o(1)}{\log^{*(r-1)} n} \right]^{d-r+1}.$$

   While a considerable amount of attention has been devoted to the study of random bootstrap percolation in the $d$-dimensional box, it has also been extensively explored on the hypercube. Let $Q_n$ be a hypercube graph whose vertices are the binary $n$-dimensional vectors and edges connect the pairs of vertices at Hamming distance $1$. Balogh and Bollob{\'a}s \cite{balogh_bollobas_2006} established the following result: there exist constants $c_1>0$ and $c_2>0$ such that
       $$ \frac{c_1}{n^2}2^{-2\sqrt{n}}\leq p_c(Q_n,2) \leq \frac{c_2}{n^2}2^{-2\sqrt{n}}.$$
    In a subsequent study, Balogh, Bollob{\'a}s, and Morris further improved upon the bound presented in \cite{balogh_bollobas_morris_2010} and proved the following result: there exists an explicit constant $\lambda > 0$ such that  
     $$ \frac{16\lambda}{n^2}\left(1+\frac{\log n }{\sqrt{n}} \right)2^{-2\sqrt{n}}\leq p_c(Q_n,2) \leq \frac{16\lambda}{n^2}\left(1+\frac{5 (\log n)^2 }{\sqrt{n}} \right)2^{-2\sqrt{n}}.$$
    For majority bootstrap percolation on the hypercube $Q_n$ where the infection threshold $r$ is equal to $\frac{1}{2}n$, Balogh, Bollob{\'a}s and Morris in \cite{balogh_bollobas_morris_2009}  proved that: 
    $$ 
    \frac{1}{2}-\frac{1}{2}\sqrt{\frac{\log n}{n}}-2\frac{\log \log n}{\sqrt{n \log n}} \leq p_c\Big(Q_n,\frac{n}{2}\Big) \leq \frac{1}{2}-\frac{1}{2}\sqrt{\frac{\log n}{n}}+\frac{\log \log n}{2\sqrt{n \log n}}+o\left(\frac{\log \log n}{2\sqrt{n \log n}} \right).
    $$   
Recently, Collares, Erde, Geisler, and Kang \cite{collares2024universalbehaviourmajoritybootstrap} showed that the random majority bootstrap percolation process on graphs satisfying certain structural properties exhibits a universal behavior, which is, in some sense, governed by the degree sequence. Subsequently, the same authors \cite{COLLARES2026104347} considered a slightly different class of graphs and established analogous results.

\emph{Overview of the results.} 
Holroyd, Liggett, and Romik \cite{holroyd} studied bootstrap percolation on a generalized grid graph in which each vertex has more than four neighbors. In Section 9 of \cite{collares2024universalbehaviourmajoritybootstrap}, the authors posed the question of determining the critical probability for majority bootstrap percolation on the halved cube, a graph on $\{0,1\}^n$ where two vertices are adjacent if their Hamming distance is exactly two. Motivated by these works, we consider an $n$-dimensional hypercube with extended neighborhoods and establish estimates for the critical probability of majority bootstrap percolation.

Let $Q_{k,n}$ be a graph with vertices $V(Q_{k,n})=\{0,1\}^n$ and edges $E(Q_{k,n})=\{xy:d_H(x,y)=1,2,\cdots,k\}$, where $k \geq 2$. We prove the following theorem.

\begin{theorem} \label{theorem:main}
    Consider a positive integer $k \geq 2$, and let $N=\sum_{i=1}^k \binom{n}{i}$. For every $\epsilon >0$, there exists an $n_0(\epsilon,k) \in \mathbb{N}$ such that for $n \geq n_0(\epsilon,k)$ we have 
    
$$  \frac{1}{2}- \left(\frac{\sqrt{kk!}}{2}+\epsilon \right)\frac{\sqrt{\log n}}{n^{k/2}} \leq p_c\left(Q_{k,n},\frac{N}{2} \right) \leq \frac{1}{2}- \left(\frac{\sqrt{kk!}}{2}-\epsilon \right)\frac{\sqrt{\log n}}{n^{k/2}}.$$

\end{theorem}
Note that our graph does not lie on the class of graphs analyzed in 
\cite{collares2024universalbehaviourmajoritybootstrap} or \cite{COLLARES2026104347} while the critical probability exhibits the same universal second-order behavior in \cite{collares2024universalbehaviourmajoritybootstrap} and in \cite{COLLARES2026104347} when expressed in terms of its degree $N$. 

\emph{Our Contributions.} Our first contribution is to apply a theorem of Petrov, a theorem of Gordon, a new lemma \ref{lemma: binom approx} for estimating the binomial distribution, and Effron-Stein's inequality to obtain very precise approximations for binomial tail probabilities. This level of precision  plays a crucial role in establishing both the upper and lower bounds for the critical probability. 

In the work by Balogh, Bollob{\'a}s, and Morris \cite{balogh_bollobas_morris_2009}, the second-moment method was used to show that there exists a $\delta >0$ such that 
\begin{align} 
    \mathbb{P}(x \in A_2) \geq \frac{1}{2}+\delta. \label{eq:1/2+delta}
\end{align}
The proof was based on the relations $\mathcal{N}(y) \cap \mathcal{N}(z)=\{x,w\}$ and $\mathcal{N}(x) \cap \mathcal{N}(y)=\emptyset$, where $y, z \in \mathcal{N}(x)$ and $x \in V(Q_n)$ which made the application of the second-moment analysis manageable. In \cite{collares2024universalbehaviourmajoritybootstrap} and \cite{COLLARES2026104347}, the authors followed the approach of \cite{balogh_bollobas_morris_2009} by using the graph properties essentially similar (in some sense) to the hypercube to derive equation \ref{eq:1/2+delta}. In our case, these relations and the graph properties are absent, so application of the second-moment analysis similar to theirs lacks a good starting point. More precisely, we have for $k \geq 2$ 
$$\max_{\substack{y,z \in \mathcal{N}(x)\\y \neq z} } |\mathcal{N}(y) \cap \mathcal{N}(z)|=\Theta(n^{k-1}),$$
and 
$$\max_{\substack{x,y\\ y \in \mathcal{N}(x)}}|\mathcal{N}(x) \cap \mathcal{N}(y)|=\Theta(n^{k-1}).$$

Let us briefly outline our proof strategy for equation~\ref{eq:1/2+delta}. Let $S_x=\sum_{v \in \mathcal{N}(x)}\mathbb{I}_{\{v \in A_1\}}$. We manage to prove $\mathbb{P}(S_x \geq \frac{N}{2})=1-o(1)$ by obtaining sufficiently accurate estimates for  $\mathbb{E}S_x$ and $\text{Var}(S_x)$. To derive the required level of precision, we apply Lemmas ~\ref{lemma: Petrov}-\ref{lemma: efron}. 

Another contribution of our paper is the application of Theorems ~\ref{theorem: Baranyai1} and ~\ref{theorem: Baranyai2} so that we can adopt the same proof strategy of ~\cite{balogh_bollobas_morris_2009} to prove 
$\mathbb{P}(x \in A_5)=1-\exp(-cn)$ where $c$ is a positive constant. More precisely, they provide a partition of $(S([0]^n,k))$, where $S(x,k)=\{y\in V(Q_{k,n}): d_H(x,y)=k\}$, such that the parts have approximately equal size and, within each part, the supports of distinct vectors are pairwise disjoint. This structure gives us effective control over the distributions of the neighbors associated with each part.

Our third contribution is a non-trivial enumeration of the possible choices for $S \in S([0]^n,4)$, which allows Lemma~\ref{lemmaA8} to be applied effectively to prove the lower bound of Theorem~\ref{theorem:main}. When $k=2$, the exponential bound provided by Lemma~\ref{lemmaA8} is not sufficiently strong if the number of possible choices for $S$ is estimated trivially. Viewing the elements of a set $G\subset S([0]^n,3))$ as the edges of a $3$-uniform hypergraph on $[n]$, Balogh, Bollob{\'a}s, and Morris \cite{balogh_bollobas_morris_2009} considered the codegree $d_G(i,j)=\bigl|{\ell: ij\ell\in E(G)}\bigr|$
of each pair $i,j\in[n]$, and showed that the presence of many pairs with large codegree substantially reduces the number of possible choices for $G$. In our setting, we refine this approach by considering not only the codegrees of pairs but also those of triples. Moreover, we establish a quantitative relationship between the codegrees of pairs and triples and the quantity $D(k)$ appearing in Lemma~\ref{lemmaA8}. Deriving this relationship is more involved than the corresponding argument in \cite{balogh_bollobas_morris_2009}. The main arguments used in \cite{collares2024universalbehaviourmajoritybootstrap} and \cite{COLLARES2026104347} to establish the lower bound on the critical probability rely heavily on structural properties of the underlying graphs that are absent in our setting. Therefore, these arguments do not appear to be directly applicable in our case.

\emph{Our Proof Strategy.} At a high level we follow the proof strategy of \cite{balogh_bollobas_morris_2009}. The main intuition is as follows. When the initial infection probability is slightly above the critical probability, after a few steps each vertex remains healthy with exponentially small probability. When the initial infection probability is slightly below the critical probability, after a few steps the infection process stops with high probability.

The heuristics behind the constant in the second-order term are as follows. Let the initial infection probability $p=\frac{1}{2}-\frac{C\sqrt{\log n}}{n^{k/2}}$. Let $\tau_n=\mathbb{P}(x \in A_1| x \notin A_0)=\mathbb{P}(|\mathcal{N}(x) \cap A_0|\geq r)$. For $r=\frac{N}{2}$, we can apply Lemmas ~\ref{lemma: Petrov} and ~\ref{lemma:Gordon} to get $\tau_n=(1+o(1))\frac{\sqrt{k!}}{2C \sqrt{2\pi \log n}} n^{\frac{-2C^2}{k!}}$. For a typical healthy vertex $x$ initially there are about $Np=\frac{N}{2}-N\frac{C\sqrt{\log n}}{n^{k/2}}$ infected neighboring vertices. The deficit to its threshold $r$ is about $N\frac{C\sqrt{\log n}}{n^{k/2}}$. Among all the neighboring vertices of $x$, there are roughly $N\tau_n$ newly infected vertices in the first step. The transition occurs when $\tau_n \approx \frac{\sqrt{\log n}}{n^{k/2}}$. 

For the upper bound of Theorem ~\ref{theorem:main}, we first show equation~\ref{eq:1/2+delta}. Then we show that after a few steps, the probability that $x$ remains healthy is exponentially small, of order $\exp(-cn)$. Next we strengthen this estimate so that we can use the union bound over all vertices. 

For the lower bound of Theorem ~\ref{theorem:main}, we introduce a less restrictive process in which for the first few steps the number of infected neighbors for a vertex $x$ to be infected is slightly less than $r$. Then we show this less restrictive process will stop in a few steps with high probability. This implies that the original process will fail to percolate with high probability.

\section{Preliminaries and Tools}
We will use the standard Chernoff bound (See Remark 2.5 in \cite{Random_graphs}) which was also used in \cite{balogh_bollobas_morris_2009}.
\begin{lemma} \label{lemmaA1}
Let $n \in \mathbb{N}$, $0<p<1$, $t>0$ and $S(n) \sim \text{Binomial}(n,p)$. Then,
$$\mathbb{P}(S(n) \geq np+t) \leq \exp\left(-\frac{2t^2}{n}\right)$$
and
$$\mathbb{P}(S(n) \leq np-t) \leq \exp\left(-\frac{2t^2}{n}\right).$$
\end{lemma}

 The following lemma is taken from \cite{balogh_bollobas_morris_2009}, which generalizes Lemma~\ref{lemmaA1} to a weighted binomial distribution. This lemma is crucial for us to prove the lower bound of Theorem \ref{theorem:main}.
\begin{lemma} \label{lemmaA8}
Let $t,k,d_1,...,d_k \in \mathbb{N}$ and $p \in (0,1)$. Let $X_i \sim \text{Binomial}(d_i,p)$ for each $i \in [k]$, let $Y_k=\sum_{i=1}^kiX_i$, and let $D(k)=\sum_{i=1}^ki^2d_i$. Then 
$$\mathbb{P}(Y_k \geq \mathbb{E}(Y_k)+t) \leq (2t)^{k-1}\exp\left(\frac{-2t^2}{D(k)} \right).$$
\end{lemma}
We will need the following fact.
\begin{fact} \label{fact3}
    If $S \subset V(Q_{k,n})$ satisfies $d_H(y,z) \geq 2k\ell+1$ for every $y,z \in S$, then the events $\{y \in A_{\ell}\}_{y \in S}$ are independent.
\end{fact}
\begin{proof}
Let $x \in V(Q_{k,n})$, define $B_G(x,\ell):=\{z\in V(Q_{k,n}): d_G(x,z)\le \ell\}.$

We prove by induction on $\ell$ that the event $\{x\in A_\ell\}$ is completely determined by the initial infection states of the vertices in $B_G(x,\ell)$.

Initial step:
The event $\{x\in A_0\}$ simply depends only on the initial state of $x$ by definition and note that $B_G(x,0)=\{x\}$. Therefore the statement is true for $\ell=0$.

Inductive step:

Suppose that, for every vertex $u$, the event ${u\in A_{\ell-1}}$ is determined by the initial states in $B_G(u,\ell-1)$.

A vertex $x$ is infected by step $\ell$ precisely when either $x\in A_0$,
or at least $r$ neighbors of $x$ belong to $A_{\ell-1}$. Therefore, to decide whether $x\in A_\ell$, we need to know, for each $u\in \mathcal{N}_G(x)$, whether $u\in A_{\ell-1}$.

By the induction hypothesis, this is determined by the initial states in $B_G(u,\ell-1)$. Hence $x\in A_\ell$ is determined by the initial states in ${x}\cup\bigcup_{u\in N_G(x)}B_G(u,\ell-1)$.

If $z\in B_G(u,\ell-1)$ and $u$ is a neighbor of $x$, then 
$$d_G(x,z)
\le d_G(x,u)+d_G(u,z)
\le 1+(\ell-1)
=\ell.$$

Thus
$$
B_G(u,\ell-1)\subseteq B_G(x,\ell),$$
so
$$
x\cup\bigcup_{u\in \mathcal{N}_G(x)}B_G(u,\ell-1)
\subseteq B_G(x,\ell).
$$

Therefore,$ {x\in A_\ell}\text{ is determined by the initial states in }B_G(x,\ell)$.

In the graph $Q_{k,n}$, one edge can change at most $k$ coordinates. In fact, $ d_G(x,z)=\left\lceil
\frac{d_H(x,z)}{k}
\right\rceil, $ and consequently $B_G(x,\ell)=\{z: d_H(x,z)\le k\ell\}$.

Hence ${x\in A_\ell}$ depends only on the initial infection variables whose Hamming distance from $x$ is at most $k\ell$.   
\end{proof}

We will also use the following three lemmas which were proved in \cite{balogh2007majority}. 
\begin{lemma} \label{lemma:S(n) m}
    Let $p \in (0,1)$ and $n \in \mathbb{N}$ satisfy $pn^2 <1$, and let $S(n) \sim \text{Bin}(n,p)$. Then 
    $$\mathbb{P}(S(n) \geq m) \leq 2p^{\frac{m}{2}}$$
for every $m \in [n]$. In particular, if $c,\epsilon_1 >0 $ and $p \leq e^{-cn}$, then for $b=b(x,\epsilon_1)$ not depending on $n$,
$$\mathbb{P}(S(n)\geq \epsilon_1 n) \leq e^{-bn^2}.$$
\end{lemma}

\begin{lemma} \label{lemma: graph partition}
Let $G$ be a graph , let $k,m \in \mathbb{N}$, and suppose that for each $x \in V(G)$,
$$B(x,k):=|y \in V(G):d(x,y) \leq k| \leq m.$$
Then there exists a partition
$$V(G)=B_1 \cup \cdots \cup B_m$$
of $V(G)$, such that if $y,z \in B_i$, then $d(y,z) \geq k+1.$ Here, $d(x,y)$ denotes the graph-theoretic distance between the vertices $x$ and $y$. 
\end{lemma}

\begin{lemma} \label{lemma: graph partition for hypercube}
 Let $n,k \in \mathbb{N}$, and $x \in V(Q_n)$. Then there exists a partition 
 $$S(x,k):=\{y \in V:d_H(y,x)=k\}=B_1 \cup \cdots \cup B_m$$
 of $S(x,k)$ into $m \leq k \binom{n}{k-1} \leq 2n^{k-1}$ sets, such that if $x,y \in B_j$ for some $j$, then $d_H(x,y) \geq 2k$. 
\end{lemma}

We need to approximate the binomial distribution very accurately. The following two lemmas will play a crucial role in our analysis. 

\begin{lemma} Chapter: Probabilities of Large Deviations in \cite{Petrov1975} \label{lemma: Petrov}

 Let $V_{N}^2=\sum_{i=1}^{N} \Var(X_{_i})$ and $S_{N}=\sum_{i=1}^{N} X_{i}$, where the centered random variables $X_{i}$  are independent with a common distribution. Assume that there exists a constant $h>0$, independent of $i$, such that 
$$\sup_{i}\mathbb{E}e^{tX_{i}}<\infty,$$
in the interval $|t| \leq h$. 
Then we have for $x \geq 0$ and $x=o(\sqrt{N})$, 
$$\mathbb{P}(S_{N} \geq x V_{N})=\left(1-\Phi(x) \right)  \exp\left(\frac{x^3}{\sqrt{N}} \lambda \left(\frac{x}{\sqrt{N}} \right) \right) \left(1+O\left(\frac{1+x}{\sqrt{N}} \right) \right),$$
where $\Phi(x)$ denotes the cumulative distribution function of a standard normal random variable,i.e.,
$$\Phi(x)=\frac{1}{2\pi}\int_{-\infty}^{x} e^{-t^2/2},$$
and 
$\lambda(t)=\sum_{k=0}^{\infty} a_kt^k$ is a power series with coefficients depending on the cumulants of the random variable $X_i$ which converges for sufficiently small values of $|t|$. 
\end{lemma}

\begin{lemma} \cite{Gordon1941} \label{lemma:Gordon}

For $x >0$, we have 
$$\frac{x}{1+x^2}\phi(x) \leq 1- \Phi(x) \leq \frac{\phi(x)}{x},$$
where $\phi(x)=\frac{1}{2\pi}e^{-x^2/2}$. 
    
\end{lemma}

Let us fix some notations. 
Let $C=\frac{\sqrt{kk!}}{2}-\epsilon$, $N=\sum_{i=1}^k \binom{n}{i}$, and $\Delta_n=\left(\frac{\sqrt{kk!}}{2}-\epsilon \right)\frac{\sqrt{\log n}}{n^{k/2}}$. 

The following lemma gives a good estimate for the probability of a binomial random variable which moderately deviates from its mean. 
\begin{lemma} \label{lemma: binom approx}
Suppose $L=N+O(n^{k-1})$ and $Y\sim \text{Binomial}(L,p)$. Suppose an integer $q$ satisfies
$$q-Lp=N\Delta_n+o(N\Delta_n).$$
Then $$
\mathbb P(Y=q) \le N^{-1/2}n^{-\alpha+o(1)},
$$
where $\alpha=\frac{2(\frac{\sqrt{kk!}}{2}-\epsilon)^2}{k!}$.

\end{lemma}
\begin{proof}
Let $\theta= \frac{q}{L}$. By Stirling’s formula we have

$$\mathbb{P}(Y=q)= \frac{1+o(1)}{\sqrt{2 \pi L \theta (1-\theta)}} \exp(-LD_{KL}(\theta || p)),$$
where $D_{KL}(\theta || p)=\theta \log \frac{\theta}{p}+(1-\theta) \log \frac{1-\theta}{1-p}$.

Let $f_p(x)=x\log \frac{x}{p}+(1-x)\log \frac{1-x}{1-p}$. Then we have 
$\frac{df_p(x)}{dx}= \log(x)-\log(p)-\log(1-x)+\log(1-p)$, $\frac{d^2f_p(x)}{d^2x}= \frac{1}{x(1-x)}$, $f(p)=0$, and $\frac{df_p(x)}{dx}|_{x=p}=0$.

Since $\theta -p=\Delta_n(1+o(1))$, the Taylor expansion around $p$ gives
$$D_{KL}(\theta || p)=\frac{(\theta -p)^2}{2p(1-p)}+O(|\theta -p|^3).$$

Since $p(1-p)=\frac{1}{4}+o(1)$,  we have$$L D_{KL}(\theta ||p)=2N\Delta_n^2+o(\log n).$$

Note that the cubic term is negligible because $N\Delta_n^3=O(n^{-k/2} (\log n)^{3/2})=o(1)$. We also have 
$2N\Delta_n^2=(\alpha+o(1)) \log n$. 

Finally note that $\frac{1}{\sqrt{2 \pi L \theta (1-\theta)}}=O(\frac{1}{N^{1/2}})$. \qedhere
\end{proof}

For our analysis, we will also need the following two lemmas. 
\begin{lemma} Efron-Stein Inequality \label{lemma: efron}

Suppose $X_1,X_2,\cdots,X_n, X_1',X_2', \cdots, X_n'$ are independent with $X_i, X_i'$ having the same distribution for all $i$. Let $X=(X_1,X_2,\cdots,X_n)$ and $X^{(i)}=(X_1,X_2,\cdots, X_{i-1},X_i',X_{i+1},\cdots,X_{n})$. Let $f(X)$ be a function of these random variables. Then 
$$\text{Var}(f(X)) \leq \frac{1}{2}\sum_{i=1}^n \mathbb{E}(f(X)-f(X^{(i)}))^2.$$
\end{lemma}

\begin{lemma} Chebyshev's inequality

Let $X$ (integrable) be a random variable with finite non-zero variance $\sigma^2$ and finite expected value $\mu$. Then for any real number $a >0$, 
$$\mathbb{P}(|x-\mu| \geq a) \leq \frac{\sigma^2}{a^2}.$$
    
\end{lemma}

Before stating the next theorem, we will introduce two definitions.
\begin{definition}
   A 1-factor of a hypergraph is a spanning 1-regular sub-hypergraph.  
\end{definition} 

\begin{definition}
 A decomposition of a hypergraph into edge-disjoint 1-factors is a 1-factorization.   
\end{definition}
The following Theorem was originally proved in \cite{Zsolt} by Baranyai.
\begin{theorem} \label{theorem: Baranyai1}
  Let $k \geq 2$. A 1-factorization of the complete k-uniform hypergraph on $n$ vertices exists if and only if $k|n$.  
\end{theorem}

We will also need the stronger version of Theorem ~\ref{theorem: Baranyai1}. 
\begin{definition}
Let $\{V_i\}_{i=1}$ be sets of vertices with $|V_i|=n_i$. By $K_{n_1,n_2,\cdots,n_\ell}^h$, we denote the hypergraph $(V,\mathscr{E})$, where $$V=\cup_{i=1} V_i,$$
$$\mathscr{E}=\{E: E \subset V,|E|=h, |V_t \cap E| \leq 1, t=1,2,\cdots, \ell \}.$$
If $n_1=n_2=\cdots=n_{\ell}=m$, the hypergraph $K_{n_1,n_2,\cdots,n_\ell}^h$ is denoted by $K_{\ell\times m}^h$.
\end{definition}

If the edge set $\{\mathscr{E}_i\}_{i=1}^s \subset \mathscr{E}$ is a partition of $\mathscr{E}$, then we say that $(V, \mathscr{E})$ is partitioned into the hypergraphs $(V,\mathscr{E}_i)$. The following Theorem was initially proved in \cite{Zsolt2} by Baranyai. 
\begin{theorem} \label{theorem: Baranyai2}
Let $a_1,a_2,\cdots,a_s$, be natural numbers such that $\sum_{j=1}^sa_j = \binom{\ell}{h}m^h$. Then the edges of $K_{\ell \times m}^h$ can be partitioned into almost regular hypergraphs $(V,\mathscr{E}_i)$ (i.e., the degree of each vertex differs by at most one) so that  $|\mathscr{E}_j|= a_j$ for $j = 1,..., s$.
\end{theorem}

\section{Upper Bound}
For the ease of notations, from now on by $c=c(k)$ we denote a constant which may depend on $k$. We will use $c=c(k)$ as generic notations, where the specific form of this function may be different in different expressions. 

Let us start by proving the upper bound. To do so, we require a slightly more general concept. Let $A_{i,r}$ represent the set of infected vertices up to and including step $i$ when the infection threshold is $r$. Hence, $A_i = A_{i,\frac{N}{2}}$. In order to establish the upper bound, we will rely on a few lemmas.

\begin{lemma} \label{lemma:0.5+0.5delta}
   Let $\epsilon >0, x \in V(Q_{k,n})$, let $p=\frac{1}{2}-\left(\frac{\sqrt{kk!}}{2}-\epsilon \right)\frac{\sqrt{\log n}}{n^{k/2}}$, and let the infection threshold be $r =\frac{N}{2}+O(1)$, where $k \geq 2$ is an integer. Then for sufficiently large $n(\epsilon)$ there exists an absolute constant $\delta >0$ such that  
    
$$\mathbb{P}(x \in A_{2,r}) \geq \frac{1}{2}+ \delta.$$ 
\end{lemma} 

In fact we will prove that under the assumptions in Lemma~\ref{lemma:0.5+0.5delta}, we have $\mathbb{P}(x \in A_{2,r}) \geq 1-o(1).$
\begin{proof} 
The main strategy is to prove that $\mathbb{P}(S_x \leq \frac{N}{2})=o(1)$, where $S_x=\sum_{u \in \mathcal{N}(x)} \mathbb{I}_{\{u \in A_{1,r}\}}$. In order to prove this statement, it is very natural to consider using the Chebyshev's inequality. Thus, we need to estimate the
$\mathbb{E}(S_x)$ and $\text{Var}(S_x)$.

First let us consider $\mathbb{P}(x \in A_{1,r})$ and so we have 
\begin{align*}
    \mathbb{P}(x \in A_{1,r}) &=\mathbb{P}(x \in A_0)+\mathbb{P}(|\mathcal{N}(x) \cap A_0| \geq r | x \notin A_0) \mathbb{P}(x \notin A_0)\\
    & = \mathbb{P}(x \in A_0)+\mathbb{P}(|\mathcal{N}(x) \cap A_0| \geq r) \mathbb{P}(x \notin A_0).
\end{align*}

Therefore $\mathbb{E}[\mathbb{I}_{x \in A_{1,r}}]=p+(1-p) \tau_n$, where $\tau_n=\mathbb{P}(|\mathcal{N}(x) \cap A_0| \geq r)$. 

Note that $|\mathcal{N}(x) \cap A_0|$ is a Binomial$(N,p)$. Let the random variable $X_N=|\mathcal{N}(x) \cap A_0|$.  Therefore, we have 
$\mu=\mu_N:=\mathbb{E}X_N=Np$ and $\sigma^2=\sigma_N^2:=p(1-p)N$. We would like to use Lemmas ~\ref{lemma: Petrov} and ~\ref{lemma:Gordon} to approximate $\tau_n$ very accurately. Before that we need to check the uniform boundedness condition. 

Note that $X_N-pN=\sum_{i=1}^{N} (Z_{i}-p)$ where $Z_i \sim \text{Bernoulli}(p)$. Define $Y_{i}=\frac{Z_{i}-p}{p(1-p)}$. It is easy to see $ 0.49 \leq p \leq \frac{3}{4}$ for sufficiently large $n$ and therefore $p(1-p) \geq \frac{3}{16}$. Since $Z_{i} \in \{0,1\}$, $|Z_{i}-p| \leq 1$. Therefore, we have 
$$|Y_{i}| \leq \frac{1}{p(1-p)} \leq \frac{4}{\sqrt{3}}.$$
Thus for every $h >0$, we have 
$$e^{h |Y_{i}|} \leq e^{\frac{4h}{\sqrt{3}}},$$
and consequently, 
$$\sup_{i}\mathbb{E}e^{h |Y_{i}|} \leq e^{\frac{4h}{\sqrt{3}}} < \infty.$$
Recall that $\Delta_n=\left(\frac{\sqrt{kk!}}{2}-\epsilon \right)\frac{\sqrt{\log n}}{n^{k/2}}$ and note that 
$\tau_n = \mathbb{P}\left(\frac{X_N-\mu_N}{\sigma_N} \geq \frac{r-\,u_N}{\sigma_N}\right)$. 

We have $r=\frac{N}{2}+O(1)$, and $pN=\frac{N}{2}-N\Delta_n$. Therefore, 
$r-\mu=N\Delta_n+O(1)$, and  $\sigma=(\frac{1}{2}+o(1))\sqrt{N}.$

Therefore by Lemmas ~\ref{lemma: Petrov} and ~\ref{lemma:Gordon}, $$\tau_n=(1+o(1))\frac{1}{x\sqrt{2\pi}}\exp(-\frac{x^2}{2}).$$
where $x:=\frac{r-\mu}{\sigma}=\left(\frac{\sqrt{kk!}-2\epsilon}{\sqrt{k!}}+o(1)\right) \log n.$

It is easy to see $\tau_n=(1+o(1))\frac{\sqrt{k!}}{2C \sqrt{2\pi \log n}} n^{\frac{-2C^2}{k!}}$ and recall that $C=\frac{\sqrt{kk!}}{2}-\epsilon$.

Note that $\tau_n=n^{-\frac{2C^2}{k!}+o(1)} \gg \frac{\sqrt{\log n}}{n^{k/2}}$. Thus $\mathbb{P}(x \in A_{1,r})=\frac{1}{2}+(\frac{1}{2}-o(1))n^{-\frac{2C^2}{k!}+o(1)}$.

\vspace{1em}
Now we need to evaluate $\text{Var}(S_x)$. This part of the proof is a bit more involved.

For every vertex $x$, let $\xi_x\sim \text{Bernoulli}(p)$ denote its initial infection state, independently. Define $I_u=\mathbb{I}_{\{u \in  A_{1,r}\}}$ and let $\alpha=\frac{2C^2}{k!}$.

Thus
$$I_u= \mathbb{I}\left\{ \xi_u=1 \ \text{or}
\sum_{x\in \mathcal{N}(u)}\xi_x\ge r \right\}.$$

Fix distinct $u,w \in V(Q_{k,n})$, and set $d=d_H(u,w)$. We will have the following claim.

\begin{claim}
$$
\text{Cov}(I_u,I_w) \le
n^{-\lceil d/2\rceil-2\alpha+o(1)} + \mathbb{I}_{{d\le k}} n^{-k/2-\alpha+o(1)}.
$$
\end{claim}

\begin{proof}
Note that the events $I_u$ and $I_v$ may have dependency due to the common neighbors i.e., $\mathcal{N}(u)\cap \mathcal{N}(w)$ and the facts that $u\in \mathcal{N}(w)$ and $w\in \mathcal{N}(u)$ when $d\le k$.

Let us estimate the Size of the common neighborhood. Let $M:=\mathcal{N}(u)\cap \mathcal{N}(w)$. We first show $M=O\left(n^{k-\lceil d/2\rceil}\right)$.

Without loss of generality, we can assume that $u=[0]^n$, and write  $W:=\{w:||w||_{\ell_1}\}=d$. For a common neighbor $x$, write $a=|x\cap W|$ and $b=|x\setminus W|$.

Then we have 
$$d_H(x,u)=a+b\le k $$
and
$$d_H(x,w)=d-a+b\le k.$$
Adding the two inequalities gives
$$d+2b\le 2k.$$
Therefore,
$$b\le k-\left\lceil\frac d2\right\rceil.$$

There are only $O(1)$ choices for $x\cap W$. Consequently,
$$
M = O_k\left( n^{k-\lceil d/2\rceil}\right).
$$
Note that $M=0$ if $d>2k$.

\textbf{Case 1:} $d >k$.

In this case, $u\notin \mathcal{N}(w)$ and $w\notin \mathcal{N}(u)$. Thus the only initial variables that influence both $I_u$ and $I_w$ are those indexed by the common neighborhood $\mathcal C$.

Let
$$ Z:=\sum_{x\in\mathcal C}\xi_x \sim \text{Binomial}(M,p).
$$

Conditional on all variables $(\xi_x)_{x\in\mathcal C}$, the remaining variables determining $I_u$ and $I_w$ are disjoint. Therefore $I_u$ and $I_w$ are conditionally independent.

Moreover, the conditional probability depends on the common variables only through their sum $Z$. Define
$$
h(s):= \mathbb{P}(I_u=1\mid Z=s).
$$

Since $u$ is initially infected with probability $p$,
$$h(s)= p+(1-p) \mathbb P\left( \text{Binomial}(N-M,p)\ge r-s
\right). $$

The same function $h$ applies to $w$. Conditional independence gives
\begin{align*}
\mathbb{E}[I_uI_w\mid Z]
&= \mathbb{E}[I_u\mid Z]\mathbb{E}[I_w\mid Z] \\
&= h(Z)^2.
\end{align*}
Thus by the tower property of expectation, we have 
$$\text{Cov}(I_u,I_w) =\text{Var}(h(Z)).$$
Therefore it remains to estimate $\text{Var}(h(Z))$.

Write
$$Z=\xi_1+\cdots+\xi_M,$$
where the $\xi_i \sim \text{Bernoulli}(p)$.

Let $\xi_i'$ be an independent copy of $\xi_i$, and let

$$Z^{(i)}= Z-\xi_i+\xi_i'.$$

The Efron–Stein inequality in Lemma~\ref{lemma: efron} gives
$$\text{Var}(h(Z)) \le\frac{1}{2} \sum_{i=1}^{M} \mathbb{E}
\bigl(h(Z)-h(Z^{(i)})\bigr)^2 $$

Let $T_i:=Z-\xi_i \sim \text{Binomial}(M-1,p)$.

If $\xi_i=\xi_i'$, then $Z=Z^{(i)}$. If they differ, then the difference in the arguments of $h$ is $1$. Therefore,

$$\mathbb{E}\bigl(h(Z)-h(Z^{(i)})\bigr)^2= 2p(1-p) \mathbb{E} \bigl(h(T_i+1)-h(T_i)\bigr)^2.
$$

Hence $$\text{Var}(h(Z)) \le Mp(1-p) \mathbb{E}\bigl(h(T+1)-h(T)\bigr)^2 $$
where $T\sim \text{Binomial}(M-1,p)$.

Note that 
$h(s+1)-h(s) =(1-p)\mathbb{P}\left(\text{Binomial}(N-M,p)=r-s-1
\right).$ Now we need to bound this quantity. 

Define the event $$\mathcal{T}= \{T: |T-(M-1)p|\le \sqrt{M} \log n
\}.$$

By Chernoff’s bound in Lemma \ref{lemmaA1}, we have 
$$\mathbb P(\mathcal{T}^c) \le e^{-\Omega(\log^2 n)}.$$

For $T\in\mathcal{T}$, let $q=r-T-1$.Then we have 

\begin{align*}
q-(N-M)p
&= r-T-1-(N-M)p\\
&= \frac{N}{2}-Np+O(1)
-\bigl(T-(M-1)p\bigr)\\
&=N\Delta_n +O(\sqrt M\log n).
\end{align*}

Since $\sqrt{M}\log n=o(N\Delta_n)$, Lemma ~\ref{lemma: binom approx} can be applied.

Therefore,
$$\mathbb P\left( \text{Binomial}(N-M,p)=q \right)
\le
N^{-1/2}n^{-\alpha+o(1)}$$
uniformly on $\mathcal{T}$.

Therefore, 
$$|h(T+1)-h(T)|
\le
N^{-1/2}n^{-\alpha+o(1)}
$$
on $\mathcal T$. On $\mathcal T^c$, we use the trivial bound $|h(T+1)-h(T)|\le1$. Hence

\begin{align*}
\mathbb E
\bigl(h(T+1)-h(T)\bigr)^2
&\le
N^{-1}n^{-2\alpha+o(1)}
+
e^{-\Omega(\log^2 n)}\\
&=
N^{-1}n^{-2\alpha+o(1)}.
\end{align*}

Consequently, 
$$\text{Var}(h(Z)) \le \frac{M}{N}n^{-2\alpha+o(1)}.$$

Since $N=\Theta(n^k)$,
$$\frac {M}{N} =O_k\left(n^{-\lceil d/2\rceil}\right).$$

Consequently, when $d>k$, 
$$
\text{Cov}(I_u,I_w)
\le
n^{-\lceil d/2\rceil-2\alpha+o(1)}.$$

This is the first term in (1).

\textbf{Case 2:} $d\leq k$

Now $u$ and $w$ are adjacent, so there is additional dependence:
$
w\in \mathcal{N}(u), u\in \mathcal{N}(w).
$

We remove this direct dependence by replacing the two cross-variables with independent copies.

Let $\widehat\xi_u,\widehat\xi_w$ be independent Bernoulli($p$) variables, independent of all original variables. Define
\begin{align*}
J_u &= 
\begin{cases}
1 & \text{if 
$\xi_u=1$ or 
$\sum_{x \in \mathcal{N}(u) \setminus w}\xi_x
+\widehat\xi_w
\ge r$} \\
0 & \text{otherwise}
\end{cases}
\end{align*}

and similarly, 

\begin{align*}
J_w &=
\begin{cases}
1 & \text{if 
$\xi_w=1$
or
$\sum_{x\in \mathcal{N}(w)\setminus{u}}\xi_x
+\widehat\xi_u
\ge r$} \\
0 & \text{otherwise}
\end{cases}
\end{align*}

Note that  $J_u\overset d=I_u,
\qquad
J_w\overset d=I_w.$

After this replacement, the only variables shared by $J_u$ and $J_w$ are the common-neighbor variables $(\xi_x)_{x\in\mathcal C}$.

Indeed, put $U=\mathcal{N}(u)\setminus(\mathcal C\cup{w})$ and $W'=\mathcal{N}(w)\setminus(\mathcal C\cup{u})$. Then $U\cap W'=\emptyset$. 

Conditional on the variables in $\mathcal C$, $J_u$ depends on $\xi_u,\widehat\xi_w$, and $(\xi_x)_{x\in U}$ and $J_w$ depends on $\xi_w,\widehat\xi_u$, and $(\xi_x)_{x\in W'}$.

These two collections are disjoint and independent. Therefore, exactly as above,
$$
\text{Cov}(J_u,J_w)
\le
n^{-\lceil d/2\rceil-2\alpha+o(1)}.
$$

It remains to compare $I_u,I_w$ with $J_u,J_w$.

The variables $I_u$ and $J_u$ can differ only when changing the neighbor bit $\xi_w$ to $\widehat\xi_w$ changes whether $u$ reaches threshold.

Let $$ R_u
:=
\sum_{x\in N(u)\setminus{w}}\xi_x
\sim\text{Binomial}(N-1,p).
$$

Then $I_u\ne J_u$ can occur only if all three of the following happen:

$$
\xi_u=0,\qquad
R_u=r-1,\qquad
\xi_w\ne\widehat\xi_w.
$$

In fact, these conditions are also sufficient. Therefore,
\begin{align*}
\mathbb P(I_u\ne J_u)
&=
(1-p)
\mathbb P\left(
\operatorname{Bin}(N-1,p)=r-1
\right)
\mathbb P(\xi_w\ne\widehat\xi_w)\\
&=
2p(1-p)^2
\mathbb P\left(
\operatorname{Binomial}(N-1,p)=r-1
\right).
\end{align*}

Now $r-1-(N-1)p= N\Delta_n+O(1).$ By Lemma ~\ref{lemma: binom approx}, we have
$$
\mathbb P\left(
\operatorname{Binomial}(N-1,p)=r-1
\right)
\le
N^{-1/2}n^{-\alpha+o(1)}.
$$

Since $N^{-1/2}=n^{-k/2+o(1)}$,$$
\mathbb P(I_u\ne J_u)
\le
n^{-k/2-\alpha+o(1)}.
$$

The same bound holds for $w$:
$$
\mathbb P(I_w\ne J_w)
\le
n^{-k/2-\alpha+o(1)}.
$$

Since $J_u=I_u$ and $J_w=I_w$, 
$$
\operatorname{Cov}(I_u,I_w)
-\operatorname{Cov}(J_u,J_w)
=
\mathbb E[I_uI_w]-\mathbb E[J_uJ_w].
$$
If $a,b,c,d \in \{0,1\}$, we have 
$$|ab-cd|
\le |a-c|+|b-d|.$$

Therefore,
\begin{align*}
\left|
\mathbb E[I_uI_w]-\mathbb E[J_uJ_w]
\right|
&\le
\mathbb E|I_u-J_u|
+
\mathbb E|I_w-J_w|\\
&=
\mathbb P(I_u\ne J_u)
+
\mathbb P(I_w\ne J_w)\\
&\le
n^{-k/2-\alpha+o(1)}.
\end{align*}

Therefore, for $d \le k$
$$
\operatorname{Cov}(I_u,I_w)
\le
n^{-\lceil d/2\rceil-2\alpha+o(1)}
+
n^{-k/2-\alpha+o(1)}. 
$$
For $d>k$, the second term is absent. Hence
$$
\operatorname{Cov}(I_u,I_w)
\le
n^{-\lceil d/2\rceil-2\alpha+o(1)}
+
\mathbb{I}_{\{d\le k\}}
n^{-k/2-\alpha+o(1)}. \qedhere
$$

\end{proof}
Recall $S_x:=\sum_{u \in \mathcal{N}(x)}I_u$ and let $d_H(u,w)=d$ for $1 \leq d \leq 2k$. Consider $\text{Var}(S_x)$ and we have 
\begin{align*}
\text{Var}(S_x) &=\sum_{u,w \in \mathcal{N}(x)}\text{Cov}(I_u,I_w)\\
& = \sum_{d=1}^{2k}N \min \left\{\binom{n}{d},N \right\}\text{Cov}(I_u,I_w)\\
&= \sum_{d=1}^{2k}O(n^{k+ \min\{d,k\}})\text{Cov}(I_u,I_w)\\
&=\sum_{d=1}^k O(n^{k+d})[n^{-\frac{d}{2}-2\alpha+o(1)}+n^{-k/2-\alpha+o(1)}]+\sum_{d=k+1}^{2k}O(n^{2k})n^{-\frac{d}{2}-2\alpha+o(1)}\\
& =O(n^{3k/2-2\alpha+o(1)})+O(n^{3k/2-\alpha+o(1)})\\
&=o((N\tau_n)^2),
\end{align*}
where $(N\tau_n)^2=\Theta(n^{2k-2\alpha +o(1)})$, and $\alpha < \frac{k}{2}$. 

Note that we have $$\mathbb{E}(S_x)=N\left(\frac{1}{2}-\frac{C\sqrt{\log n}}{n^{k/2}} \right)+N\left(\frac{1}{2}+\frac{C\sqrt{\log n}}{n^{k/2}} \right)\tau_n$$
and thus 
\begin{align*}
    \mathbb{P}(S_x \leq \frac{N}{2}) &\leq \mathbb{P}(|S_x-\mathbb{E}S_x| \geq \frac{1}{3}N\tau_n)\\
    & \leq \frac{\text{Var}(S_x)}{(\frac{1}{3}N\tau_n)^2}\\
    &= o(1). 
\end{align*}
Therefore, 
$\mathbb{P}(x \in A_2)=1-o(1)$. 
\end{proof}

Now we are ready to state the next lemma.
\begin{lemma} \label{lemma: 1-exp(-cn)}
   Let $\delta > 0$ be a real number and $k \geq 2$ be an integer. Assume that with the initial infection probability $p=\frac{1}{2}-\left(\frac{\sqrt{kk!}}{2}-\epsilon \right)\frac{\sqrt{\log n}}{n^{k/2}}$ we can establish $\mathbb{P}\left(x \in A_{2,\frac{N}{2}+\binom{3k}{2k}} \right) \geq \frac{1}{2}+\delta$. Then there exists a constant $c > 0$  such that 
   $$\mathbb{P}(x \in A_{5,\frac{N}{2}}) \geq 1-e^{-cn}.$$

\end{lemma} 

\begin{proof} Let $x \in V(Q_{k,n})$. First, assume that $k$ divides $n$. Then there exists a partition $S(x,k)=\cup_{i=1}^{m}B_i$ such that $\text{supp}(x) \cap \text{supp}(y) = \emptyset$ for all $x, y \in B_i$ and all $i \in [m]$. Moreover, $|B_i|=\frac{n}{k}$ for all $i \in [m]$ and $m=\binom{n-1}{k-1}$.

Indeed, we can consider $S(x,k)$ as the edges of a $k$-uniform hypergraph on $n$ vertices. According to Theorem~\ref{theorem: Baranyai1}, such a partition of $S(x,k)$ exists given that $k|n$.

Now assume that $k$ does not divide $n$. By Theorem~\ref{theorem: Baranyai2}, given positive integers $b_1,b_2,\cdots,b_m$ such that $\sum_{i=1}^m a_i=\binom{n}{k}$, the family $\binom{[n]}{k}$ can be partitioned into 
$$B_1 \cup B_2 \cdots \cup B_m,$$
with $|B_i|=b_i$ where each $B_i$, regarded as a $k$-uniform hypergraph, is almost regular. We can choose the $b_i$'s so that $b_i \in \{q,q-1\}$. 

Note that every $B_i$ is a matching. Indeed,
$$\frac{k|B_i|}{n} \leq \frac{kq}{n} <1,$$
since $k$ does not divide $n$. Because $B_i$ is almost regular, all its vertex degrees differ by at most $1$. Since the average degree is less than $1$, every vertex has degree $0$ or $1$. Therefore each $B_i$ consists of pairwise disjoint $k$-sets. Thus we have for each $i$
$$\left||B_i|-\frac{n}{k} \right|\leq 1.$$

From Lemma~\ref{lemmaA1} we have $$\mathbb{P}(|B_i \cap A_0| \leq \frac{n}{2k}-n^{\frac{1.1}{2}}) \leq \exp(-cn^{0.1}).$$

Let $J(i,x)$ denote this event does not happen, i.e., $|B_i \cap A_0|\geq \frac{n}{2k}-n^{\frac{1.1}{2}}$. From now on, we assume that $J(i,x)$ holds. Note that this event will help us to estimate the size of the set $S$ (defined below) which is crucial for our proof.

There exists an $i \in [m]$ such that there is a set $S \subset B_i \backslash A_0$ of size $|S|=\frac{n}{2k}+1$. If such a set $S$ cannot be found, this implies that for every $i \in [m]$, there does not exist an $S$ such that $S \subset B_i \setminus A_0$ and $|S|=\frac{n}{2k}+1$. That means for every $i \in [m]$ every ($\frac{n}{2k}$+1)-element subset of $B_i$ contains at least one element of $A_0$. Therefore for every $i \in [m]$, $|B_i \cap A_0| \geq \frac{n}{2k}+1$, then $x \in A_1$, and we are done. Indeed, since each support has size $k$, $kb_i \leq n$. Thus 
$$\frac{1}{2}\binom{n}{k}=\frac{1}{2}\sum_{i=1}^mb_i \leq m\frac{n}{2k}$$

Since the ``+1'' term will not affect the calculations, we can safely assume that there exists a set $S \subset B_i \backslash A_0$ with $|S|=\frac{n}{2k}$. Then we will show that $S \cap A_4 = \emptyset$ is highly unlikely.

Let $T(k)=\mathcal{N}(S) \cap S(x,2k)$. Denote $T_i(k)=\{y \in T(k): |\mathcal{N}(y) \cap S|=i\}$ for $i\in \{1,2\}$. It is easy to see that for $k \geq 4$
$$|T_1(k)|=\binom{n-k}{k}\frac{n}{2k}-2\binom{\frac{n}{2k}}{2}=\frac{n^{k+1}}{2k k!}-\frac{(3k-1)n^k}{4k!}+\frac{(k-1)(27k^2-19k+2)}{48k!}n^{k-1}+O(n^{k-2})$$ and 
$$|T_2(k)|=\binom{\frac{n}{2k}}{2}=\frac{n^2}{8k^2}-\frac{n}{4k}.$$
For $k=2,3,$ we need more precise calculations for $T_1(k)$. We have 
$$|T_1(2)|=\frac{n^3}{8}-\frac{11}{16}n^2+n,$$
and 
$$|T_1(3)|=\frac{n^4-12n^3+46n^2-54n}{36},$$

Let $|T_1(k) \cap A_3|=a(k)$ and $|T_2(k) \cap A_3|=b(k)$. Note that if $S \cap A_4=\emptyset$, then for $k \geq 3$ $$a(k)+2b(k) < \frac{n}{2k}\left(\frac{1}{2}N \right)=\frac{n^{k+1}}{4k k!}+\frac{(3-k)n^k}{8k(k-1)!}+\frac{3k^2-19k+50}{96k(k-2)!}n^{k-1}+O(n^{k-2})$$ and for $k=2$,
$$a(2)+2b(2) < \frac{n^3}{16}+\frac{n^2}{16}$$

We will separate the cases where $k \geq 4$, $k=3$, and $k=2$ for the clarity even though they are almost the same. 

\textbf{Case 1}: $k \geq 4$. 

$$\mathbb{E}|T_1 \cap A_0| =\frac{n^{k+1}}{4k k!}-\frac{(3k-1)n^k}{8k!}+O(n^{k-1}\sqrt{\log n}).$$

By Lemma~\ref{lemmaA1}, that there exits $c >0$ such that 
$$\mathbb{P}(||T_1 \cap A_0|-\frac{n^{k+1}}{4k k!}| \geq kn^k) \leq \exp(-cn^{k-1}).$$
Let $G(S)$ denote this event doesn't happen, i.e., $||T_1 \cap A_0|-\frac{n^{k+1}}{4k k!}| \leq kn^k $ and we have $\mathbb{P}(G(S)) \leq 1-\exp(-c n^{k-1})$.

From now on, we assume that the event $G(S)$ holds. If $|T_1\cap A_3 \backslash A_0| \geq 2kn^k$, then $a \geq \frac{n^{k+1}}{4kk!}+kn^k$ and thus $S \cap A_4 \neq \emptyset$.

\begin{claim} If $G(S)$ holds, then 
$$\mathbb{P}(|T_1\cap A_3 \backslash A_0| =O(n^k)) \leq \exp(-c\delta^2 \epsilon' n).$$
\end{claim}
\begin{proof} Consider a bipartite graph $H$ with $V(H)=H_1 \cup H_2$ where $H_1= T_1 \backslash A_0$, $H_2=S(x,3k)$ and $E(H)=\{uv: uv \in E(Q_{k,n})\}$.We will color edges of $H$ red and blue. If the end point in $H_2$ is in $A_2$, then color edges $xy$ red for every $x \in H_1$ and if the end point in $H_2$ is not in $A_2$, then color edges $xy$ blue for every $x \in H_1$.

It is easy to see that 
$|E(H)|=\binom{n-2k}{k}|H_1|$. Since $G(s)$ holds, we have $|H_1|=\frac{n^{k+1}}{4kk!}+O(n^k)$.

Now suppose that $|T_1 \cap A_3 \backslash A_0|=O(n^k)$. Since at most $O(n^k)$ vertices in $H_1$ can have at least $\frac{N}{2}$ neighbors in $A_2$, the number of red edges in $H$, denoted $E_R(H)$, satisfies 
$$|E_R(H)| \leq \frac{N}{2}|H_1|+O(n^{2k}),$$
and thus the number of blue edges in $H$, denoted by $E_B(H)$, has to satisfy
\begin{align} \label{eq: blue edges k}
    |E_B(H)| \geq |H_1|\binom{n-2k}{k}-\frac{N}{2}|H_1|-O(n^{2k}). 
\end{align}

From Lemmas \ref{lemma: graph partition} and \ref{lemma: graph partition for hypercube}, we can observe that there exists a partition of the set $S(x,3k)$ as follows:
$$S(x,3k)= \cup_{i=1}^m D_i,$$
where, for each $i \in [m]$, and for any $y$ and $z \in D_i$, we have $d_H(y,z) \geq 4k+1$. Furthermore, the number of subsets $m$ in this partition satisfies 
$$m \leq \sum_{a,b: a+b=3k \; \text{and} \; 3k-a+b \leq 4k}\binom{3k}{3k-a}\binom{n-3k}{b} \leq c n^{2k}.$$

Let us define, for $i \in [\binom{3k}{2k}]$ and $j \in [m]$,
$$D_j(i)=\{x \in D_j: |\mathcal{N}(x) \cap T_1|=i  \}$$

\begin{claim} If $ \epsilon' >0$ is arbitrarily small, then there exists an $i \in [\binom{3k}{2k}]$ and a $j \in [m]$ such that 

(a) $l$ edges of $H$ are incident with $D_j(i)$ where $l \geq \epsilon' n$ and 

(b) at most $(\frac{1}{2}+\frac{\delta}{2})l$ of those edges are red.
\end{claim}
\begin{proof} Suppose the claim is false. Let us count the number of blue edges. We have 
$$|E_B(H)| \leq \binom{n-2k}{k}|H_1|(\frac{1}{2}-\frac{\delta}{2})+\binom{3k}{2k}C(k)n^{2k} \epsilon' n,$$
which is contradiction to (\ref{eq: blue edges k}) since $\epsilon'$ is chosen arbitrarily small and $n \gg k$. 
\end{proof}
Note that events $\{ y \in A_2\}_{y \in D_j(i)}$ are independent, but they are not independent of which members of $B(x,2k)$ are in $A_0$. Since $y \in D_j(i)$ has at most $\binom{3k}{2k}$ neighbors in $B(x,2k)$ and $\mathbb{P}(y \in A_{2,\frac{N}{2}+\binom{3k}{2k}}) \geq \frac{1}{2}+\delta $, then $\mathbb{P}(y \in A_{2,\frac{N}{2}}) \geq \frac{1}{2}+\delta$ for any set $B(x,2k) \cap A_0$.

Let $E_j(i)$ denote the event that $D_j(i)$ satisfies (a) and (b). Thus 
\begin{align*}
   \mathbb{P}(\cup_{i,j}E_j(i)) & \leq \sum_{i,j}\mathbb{P}\left(\text{Bin}\left(\frac{\epsilon' n}{i},\frac{1}{2}+\delta \right) \leq \left(\frac{1}{2}+\frac{\delta}{2} \right) \frac{\epsilon' n}{i} \right) \\
   & \leq \binom{3k}{2k}n^{2k} \exp(-c\delta^2 \epsilon' n),
\end{align*}
for $c >0$.
\end{proof}

The number of choices for $S$ is at most
\begin{align*}
    \binom{\frac{n}{2k}+n^{0.55}}{\frac{n}{2k}} &\leq \left(e \left(1+\frac{n^{0.45}}{2k} \right) \right)^{n^{0.55}} \leq \exp(n^{0.6}).
\end{align*}
Let $M(S)$ denote the event $|T_1\cap (A_3 \backslash A_0)| \geq 2kn^k$. Note that if for every $i \in [m]$ and every subset $S \subset B_i \backslash A_0$ of size $\frac{n}{2k}$ all three events $J(i,x)$, $G(S)$ and $M(S)$ hold, then $x \in A_5$. 

Therefore,
\begin{align*}
   \mathbb{P}(x \notin A_5) &\leq \sum_{i}\mathbb{P}(J(i,x)^c)+\sum_{i}\sum_{S}\mathbb{P}[G(S)^c]+\sum_{i}\sum_{S}\mathbb{P}[M(S)^c] \leq \exp(-cn),
\end{align*}
where $c>0$. 

\textbf{Case 2}: $k=3$
$$\mathbb{E}|T_1 \cap A_0|=\frac{n^4}{72}-\frac{n^3}{6}+ O(n^{5/2} \sqrt{\log n}). $$
From Lemma~\ref{lemmaA1}, we have 
$$\mathbb{P}(||T_1 \cap A_0|-\frac{n^4}{72}| \geq 2n^3 ) \leq \exp(-cn^2).$$
Let $G(S)$ denote the event that $||T_1 \cap A_0|-\frac{n^4}{72}| \leq 2n^3$ and so $\mathbb{P}(G(S)) \geq 1-\exp(-cn^2)$. Assume that $G(S)$ holds. Therefore, if $|T_1 \cap A_3 \backslash A_0| \geq 3n^3$, then with high probability, $a(3) \geq \frac{n^4}{72}+n^3$ and thus $S\cap A_4 \neq \emptyset.$

\begin{claim} If $G(S)$ holds, then 
$$\mathbb{P}(|T_1 \cap A_3 \backslash A_0|=O(n^3)) \leq \exp(-c\delta^2\epsilon n).$$
\end{claim}
\begin{proof} 
Consider a bipartite graph $H$ with the vertex set $H_1 \cup H_2$, where $H_1=T_1 \backslash A_0 \subset S(x,6)$ and $H_2=S(x,9)$, and the edge set $\{uv: uv \in E(Q_{3,n}) \}$. We have 
$$|E(H)|=|H_1|\binom{n-6}{3}$$ and 
$$|H_1|=\frac{n^4}{72}+O(n^3),$$ since $G(S)$ holds.

We will color the edges of $H$ red and blue. If $y \in H_2$ is in $A_2$, then color $xy$ red for every $x \in H_1$ and if $y \in H_2$ is not in $A_2$, then color $xy$ blue for every $x \in H_1$. Now suppose $|T_1 \cap A_3 \backslash A_0|=O(n^3)$. Thus, only $O(n^3)$ vertices of $H_1$ can have at least $N/2$ neighbors in $A_2$. Thus 
$$|E_R(H)| \leq \frac{N}{2}|H_1|+O(n^6)$$
and 
\begin{align} \label{eq:blue edges 3}
    |E_B(H)| \geq \binom{n-6}{3}|H_1|-\frac{N}{2}|H_1|-O(n^6),
\end{align}
where $E_R(H)$ denotes the red edges of $H$ and $E_B(H)$ denotes the blue edges of $H$.

From Lemmas \ref{lemma: graph partition} and \ref{lemma: graph partition for hypercube}, we can observe that there exists a partition of the set $S(x,9)$ as follows:
$$S(x,9)= \cup_{i=1}^m D_i,$$
where, for each $i \in [m]$, and for any $x$ and $y \in B_i$, we have $d_H(x,y) \geq 13$. Furthermore, the number of subsets $m$ in this partition satisfies $m \leq cn^6$.

Define a refinement of $D_i$ as $D_j(i)=\{z \in D_j: |\mathcal{N}(z) \cap T_1|=i\}$ for $i \in [84].$

\begin{claim} For a sufficiently small $\epsilon''' >0$, there exist an $i \in[84]$ and $j \in [m]$ such that:

(a) $l$ edges of $H$ are incident to $D_j(i)$ and $l \geq \epsilon''' n$. And

(b) at most $(\frac{1}{2}+\frac{\delta}{2})l$ of those edges are red.
\end{claim}
\begin{proof}
Suppose the claim is false. We can have at most 
$$|E_B(H)| \leq |E(H)|(\frac{1}{2}-\frac{\delta}{2})+84\epsilon''' nm$$ blue edges, which is a contradiction to the bound in (\ref{eq:blue edges 3}) since $\epsilon'''$ is sufficiently small and $n$ is large.
\end{proof}

Note that $\{ y \in A_2\}_{y \in D_j(i)}$ are independent, but this event is not independent of which members of $B(x,6)$ are in $A_0$. Since $y \in D_j(i)$ has 84 neighbors in $B(x,6)$ and $\mathbb{P}(x \in A_{2,r=\frac{N}{2}+84}) \geq \frac{1}{2}+\delta$, then $\mathbb{P}(y \in A_{2,\frac{N}{2}}) \geq \frac{1}{2}+\delta$ for all set $B(x,6) \cap A_0$.

Let $E_j(i)$ denote the event that $D_j(i)$ satisfies (a) and (b). Thus, there exists a constant $c>0$ such that
\begin{align*}
   \mathbb{P}(\cup_{i,j}E_j(i)) & \leq \sum_{i,j}\mathbb{P}\left(\text{Bin}\left(\frac{\epsilon n}{i},\frac{1}{2}+\delta \right) \leq \left(\frac{1}{2}+\frac{\delta}{2} \right)\frac{\epsilon n}{i} \right) \\
   & \leq  \exp(-c\delta^2 \epsilon n). \qedhere
\end{align*}

\end{proof}

The number of choices for $S$ is at most
\begin{align*}
    \binom{\frac{n}{6}+n^{0.55}}{\frac{n}{6}} &\leq \exp(n^{0.6})
\end{align*}
Let $M(S)$ denote the event $|T_1\cap (A_3 \backslash A_0)| \geq 2n^3$. Note that if for every $i \in [m]$ and for every subset $S \subset B_i \backslash A_i$ of size $\frac{n}{6}$, $J(i,x)$,$G(S)$, and $M(S)$ hold, then $x \in A_5$. 

Therefore,
\begin{align*}
   \mathbb{P}(x \notin A_5) &\leq \sum_{i}\mathbb{P}(J(i,x)^c) +\sum_{i}\sum_{S}\mathbb{P}[G(S)^c]+\sum_{i}\sum_{S}\mathbb{P}[M(S)^c]\\
      & \leq \exp(-cn). \qedhere
\end{align*}

\textbf{Case 3}: $k=2$
$$\mathbb{E}|T_1 \cap A_0|=\frac{n^3}{16}-(1-\epsilon)n^2 \sqrt{\log n}+ O(n \sqrt{\log n}). $$
From Lemma~\ref{lemmaA1}, we have 
$$\mathbb{P}(||T_1 \cap A_0|-\frac{n^3}{16}| \geq 2n^2 \sqrt{\log n}) \leq \exp(-cn \log n).$$
Let $G(S)$ denote the event that $||T_1 \cap A_0|-\frac{n^3}{16}| \leq 2n^2 \sqrt{\log n}$ and so $\mathbb{P}(G(S)) \geq 1-\exp(-cn \log n)$. Assume that $G(S)$ holds. Therefore, if $|T_1 \cap A_3 \backslash A_0| \geq 3n^2 \sqrt{\log n}$, then with high probability, $a(2) \geq \frac{n^3}{16}+n^2 \sqrt{\log n}$ and thus $S\cap A_4 \neq \emptyset.$

\begin{claim} If $G(S)$ holds, then 
$$\mathbb{P}(|T_1 \cap A_3 \backslash A_0|=O(n^2 \sqrt{\log n})) \leq \exp(-c\delta^2\epsilon n).$$
\end{claim}
\begin{proof} 
Consider a bipartite graph $H$ with the vertex set $H_1 \cup H_2$, where $H_1=T_1 \backslash A_0 \subset S(x,4)$ and $H_2=S(x,6)$, and the edge set $\{uv: uv \in E(Q_{2,n}) \}$. We have 
$$|E(H)|=|H_1|\binom{n-4}{2}$$ and 
$$|H_1|=\frac{n^3}{16}+O(n^2 \sqrt{\log n}),$$ since $G(S)$ holds.

We will color the edges of $H$ red and blue. If $y \in H_2$ is in $A_2$, then color $xy$ red for every $x \in H_1$ and if $y \in H_2$ is not in $A_2$, then color $xy$ blue for every $x \in H_1$. Now suppose $|T_1 \cap A_3 \backslash A_0|=O(n^2\sqrt{\log n})$. Thus, only $O(n^2 \sqrt{\log n})$ vertices of $H_1$ can have at least $\frac{n^2+n}{4}$ neighbors in $A_2$. Thus 
$$|E_R(H)| \leq \frac{n^2}{4}|H_1|+O(n^4 \sqrt{\log n})=\frac{n^5}{64}+O(n^4 \sqrt{\log n})$$
and 
\begin{align} \label{eq:blue edges 2}
    |E_B(H)| \geq \binom{n-4}{2}|H_1|-\frac{n^2}{4}|H_1|-O(n^4 \sqrt{\log n})=\frac{n^5}{64}-O(n^4 \sqrt{\log n}),
\end{align}
where $E_R(H)$ denotes the red edges of $H$ and $E_B(H)$ denotes the blue edges of $H$.

From Lemmas \ref{lemma: graph partition} and \ref{lemma: graph partition for hypercube}, we can observe that there exists a partition of the set $S(x,6)$ as follows:
$$S(x,6)= \cup_{i=1}^m D_i,$$
where, for each $i \in [m]$, and for any $x$ and $y \in B_i$, we have $d_H(x,y) \geq 9$. Furthermore, the number of subsets $m$ in this partition satisfies $m=\binom{6}{2} \binom{n-6}{4}+\binom{6}{4} \binom{n-6}{2}+\binom{6}{3} \binom{n-6}{3}+\binom{6}{5} \binom{n-6}{1}+1$.

Define a refinement of $D_i$ as $D_j(i)=\{z \in D_j: |\mathcal{N}(z) \cap T_1|=i\}$ for $i \in [15].$

\begin{claim} For a sufficiently small $\epsilon'' >0$, there exist an $i \in[15]$ and $j \in [m]$ such that:

(a) $l$ edges of $H$ are incident to $D_j(i)$ and $l \geq \epsilon'' n$. And

(b) at most $(\frac{1}{2}+\frac{\delta}{2})l$ of those edges are red.
\end{claim}
\begin{proof}
Suppose the claim is false. We can have at most 
$$|E_B(H)| \leq |E(H)|(\frac{1}{2}-\frac{\delta}{2})+15\epsilon'' nm$$ blue edges, which is a contradiction to the bound in (\ref{eq:blue edges 2}) since $\epsilon''$ is sufficiently small and $n$ is large.
\end{proof}

Note that $\{ y \in A_2\}_{y \in D_j(i)}$ are independent, but this event is not independent of which members of $B(x,4)$ are in $A_0$. Since $y \in D_j(i)$ has 21 neighbors in $B(x,4)$ and $\mathbb{P}(x \in A_{2,r=\frac{N}{2}+21}) \geq \frac{1}{2}+\delta$, then $\mathbb{P}(y \in A_{2,\frac{N}{2}}) \geq \frac{1}{2}+\delta$ for all set $B(x,4) \cap A_0$.

Let $E_j(i)$ denote the event that $D_j(i)$ satisfies (a) and (b). Thus, there exists a constant $c >0$ such that
\begin{align*}
   \mathbb{P}(\cup_{i,j}E_j(i)) & \leq \sum_{i,j}\mathbb{P}\left(\text{Bin}\left(\frac{\epsilon n}{i},\frac{1}{2}+\delta \right) \leq \left(\frac{1}{2}+\frac{\delta}{2} \right)\frac{\epsilon n}{i} \right) \\
   & \leq 15n^4 \exp(-c\delta^2 \epsilon n). \qedhere
\end{align*}

\end{proof}

The number of choices for $S$ is at most
\begin{align*}
    \binom{\frac{n}{4}+n^{0.55}}{\frac{n}{4}} &\leq \exp(n^{0.6})
\end{align*}
Let $M(S)$ denote the event $|T_1\cap (A_3 \backslash A_0)| \geq 2n^2$. Note that if for every $i \in [m]$ and for every subset $S \subset B_i \backslash A_i$ of size $\frac{n}{4}$, $J(i,x)$,$G(S)$, and $M(S)$ hold, then $x \in A_5$. 

Therefore,
\begin{align*}
   \mathbb{P}(x \notin A_5) &\leq \sum_{i}\mathbb{P}(J(i,x)^c) +\sum_{i}\sum_{S}\mathbb{P}[G(S)^c]+\sum_{i}\sum_{S}\mathbb{P}[M(S)^c]\\
   & \leq \exp(-cn). \qedhere
\end{align*}
\end{proof}
In order to use the union bound over all vertices, We need the next lemma to strengthen the estimate for $\mathbb{P}(x \in A_\ell)$ to exponentially small error probability with exponent $\gg n$. 

\begin{lemma} \label{lemma:e^(-n^{1+k})}
   Suppose that for every $x \in V(Q_{k,n})$ there exists an absolute constant $c >0$ such that 
$$\mathbb{P}(x \notin A_{\ell}) \leq e^{-cn}.$$
Then for some $d >0$, we have 
$$\mathbb{P}(x \notin A_{(k+1)\ell+k}) \leq e^{-dn^{1+k}}.$$ 
\end{lemma} 

\begin{proof} Suppose that $x \notin A_{(k+1)\ell+k}$. We claim that for each $t \in [0,k\ell+k]$ and $t$ divisible by $k$, there exists a set $T(t) \subset S(x,t)$ and a function $ 0< \alpha(t) < 1$ such that 
\begin{align} \label{eq: T(t)}
T(t) \cap A_{(k+1)\ell+k-t}=\emptyset
\end{align}
 and 
 \begin{align} \label{eq: T(t)2}
     |T(t)|\geq \alpha(t) n^t-O(n^{t-1}).
 \end{align}

We will prove the claim by induction on $t$. When $t=0$, $|T(t)|=1$ so the statement is trivially true. Suppose that $T(t)$ satisfies \ref{eq: T(t)} and \ref{eq: T(t)2}. Then for every $y \in T(t)$ it has at most $\frac{N}{2}$ neighbors in $S(x,t+k) \cap A_{(k+1)\ell-t}$ and so at least $\binom{n-t}{k}-\frac{N}{2}$ neighbors in $S(x,t+k) \backslash A_{(k+1)\ell-t}$. Note that each element of $S(x,t+k)$ has $\binom{t+k}{k}$ neighbors in $S(x,t)$.

Therefore, there exists a set $T(t+k) \subset S(x,t+k)$ such that $$T(t+k) \cap A_{(k+1)\ell-t}=\emptyset$$ and 
$$|T(t+k)| \geq \frac{[\binom{n-t}{k}-\frac{N}{2}] |T(t)|}{\binom{t+k}{k}}.$$

From Fact~\ref{fact3}, it is easy to see that if $S \subset V(Q_{k,n})$ satisfies $d_H(y,z) \geq 2k\ell+1$ for every $y,z \in S$, then the events $\{y \in A_{\ell}\}_{y \in S}$ are independent.

From Lemmas \ref{lemma: graph partition} and \ref{lemma: graph partition for hypercube}, we can observe that there exists a partition of the set $S(x,k\ell+k)$ as follows:
$$S(x,k\ell+k)= \cup_{i=1}^m B_i,$$
where, for each $i \in [m]$, and for any $x$ and $y \in B_i$, we have $d_H(x,y) \geq 2k\ell+1$. Furthermore, the number of subsets $m$ in this partition satisfies $m =O(n^{kl})$. Indeed consider an element $y \in S(x,k\ell+k)$ and assume that $x=[0]^n$. Consider another element $z \in S(x,k\ell+k)$ and suppose $|\text{support}(y) \cap \text{support}(z)|=a$ and $|\text{support}(y) \setminus \text{support}(z)|=b$. Then we have 
$$a+b=k\ell+k$$ and 
$$k\ell+k-a+b \leq 2k\ell.$$
Therefore it is enough to take $m= |B(y,2k\ell) \cap S(x,k\ell+k)|$ and thus $m= O(n^{k\ell})$.

Now we claim that there exist a $j \in [m]$ and a sufficiently small $\epsilon_2 >0$ such that 
$$|B_j| \geq \epsilon_2 n^k$$ and 
$$|T(k\ell+k)| \geq \epsilon_2|B_j|.$$
Indeed, if the claim is not true then 
$$|T(k\ell+k)| \leq \epsilon_2 mn^k+ \epsilon_2 \binom{n}{k\ell+k},$$ which is a contraction since $\epsilon_2 $ is sufficiently small.

Recall that  $T(k\ell+k) \cap A_{\ell} =\emptyset$ and  $\mathbb{P}(x \notin A_{\ell}) \leq e^{-cn}$. From Lemma 7 in \cite{balogh2007majority}, we have
\begin{align*}
   \mathbb{P}(|T(k\ell+k) \cap B_j| \geq \epsilon|B_j|) & \leq 2(e^{-cn})^{\epsilon|B_j|/2}\\
   &\leq 2(e^{-cn^{1+k}\epsilon^2 /2}).
\end{align*}
Therefore, 
\begin{align*}
   \mathbb{P}(x \notin A_{(k+1)\ell+k}) &\leq \mathbb{P}\left( \exists j \; \text{with} \; |B_j|\geq \epsilon_2 n^k \; \text{and} \; |T(k\ell+k) \cap B_j|\geq \epsilon_2|B_j|\right)\\
   & \leq m2(e^{-cn^{1+k}\epsilon^2 /2})\\
   & \leq e^{-dn^{1+k}}. \qedhere
\end{align*}
\end{proof}
Now we are ready to prove the upper bound of Theorem~\ref{theorem:main}. 
\begin{proof}

By Lemmas ~\ref{lemma:0.5+0.5delta}, ~\ref{lemma: 1-exp(-cn)}, and ~\ref{lemma:e^(-n^{1+k})} we have 
\begin{align*}
    \mathbb{P}(A_0 \; \text{does not percolate} ) & \leq \mathbb{P}(\exists x: x \notin A_{6k+5}))\\
    & \leq \sum_{x \in V(Q_{k,n})} \mathbb{P}(x \notin A_{6k+5})\\
    & \leq 2^ne^{-dn^{1+k}}\\
    & =o(1). \qedhere
\end{align*} 
\end{proof}

\section{Lower Bound}
In this section we will prove the lower bound in Theorem~\ref{theorem:main}. Following the ideas from \cite{balogh_bollobas_morris_2009} let us introduce a new process closely related to the bootstrap percolation process. The only difference is a less strict requirement for a vertex to be infected at step 1 and step 2. Let $A_i$ denote the set of infected vertices up to step $i$. We define the process Boot($t$) as follows:

* Initially, each vertex is infected with probability $p$, independent of other vertices.

* If a vertex becomes infected, then it remains infected forever.

*A healthy vertex will be infected at step 1 if it has at least $\frac{N}{2}-2t$ infected neighbors at the initial step, i.e., $x \in A_1$ if $x \in A_0$ or $|\mathcal{N}(x) \cap A_0| \geq \frac{N}{2}-2t$.

*A healthy vertex will be infected at step 2 if it has at least $\frac{N}{2}-t$ infected neighbors at the first step, i.e., $x \in A_2$ if $x \in A_1$ or $|\mathcal{N}(x) \cap A_1| \geq \frac{N}{2}-t$.

* Let $i \geq 2$. A healthy vertex will be infected at step $i+1$ if it has at least $\frac{N}{2}$ infected neighbors at step $i$, i.e., $x \in A_{i+1}$ if $x \in A_i$ or $|\mathcal{N}(x) \cap A_i| \geq \frac{N}{2}$.

Let us refer to the original $\frac{N}{2}$-neighbor process as Boot($0$).

To establish the lower bound, we will use Lemma~\ref{lemmaA8} to show that $\mathbb{P}(A_3=A_2)=1-o(1)$. Subsequently, we will establish that $\mathbb{P}(A_2=V)=o(1).$

To prove that $\mathbb{P}(A_3 = A_2) = 1 - o(1)$, we consider separately the case $k = 2$ and the case $k \geq 3$. When $k\geq 3$, trivial bounds on the number of choices of all sets appearing in the union bound are sufficient. In contrast, when $k=2$, the trivial bound for the number of choices of the set $S \in S(x,4)$ is too weak for the union bound to succeed, and a more refined analysis is required. We first treat the case $k\geq 3$.

\begin{lemma} \label{lemmaA17}
 Consider the Boot($t$) process, where the initial infection probability $p= \frac{1}{2}-(\frac{\sqrt{kk!}}{2}+\epsilon) \frac{ \sqrt{\log n}}{n^{k/2}}$.  Let $k \geq3$ and $t=cn^{k/2}$ where $c$ is a large constant only depending on $k$. Then $A_3=A_2$ with high probability.
\end{lemma}

\begin{proof}
  Suppose $x \in A_3 \backslash A_2$. Then we have 
$$|\mathcal{N}(x) \cap A_2| \geq \frac{N}{2}$$ and 
$$|\mathcal{N}(x) \cap A_1| \leq \frac{N}{2}-t.$$
  
Thus there exists a set $T \subset \mathcal{N}(x)$ with $|T|=t$ and an $\ell \in \{ \lceil \frac{k}{2} \rceil,\lceil\frac{k}{2} \rceil+1,\cdots,k \}$ such that $T \subset A_2 \backslash A_1$ and 
\begin{equation} \label{eq: T cap S}
|T \cap S(x,\ell)| \geq \frac{t}{2k}.
\end{equation}
Indeed $|S(x,\ell)|=\binom{n}{\ell}$ for $\ell \leq n$ so as long as $n$ is sufficiently large (\ref{eq: T cap S}) is valid. 

Let $y \in T$ and then we have 
$$|\mathcal{N}(y) \cap A_1| \geq \frac{N}{2}-t$$ and 
$$|\mathcal{N}(y) \cap A_0| \leq \frac{N}{2}-2t.$$
Thus every $y \in T$ must have at least $t$ neighbors in $A_1 \backslash A_0$. Specifically, for every $y \in (T \cap S(x,\ell)),$ it must have at least $t$ neighbors in $A_1 \backslash A_0$ among which there exists an $i \in \{1,2,\cdots,k \}$ such that at least $\frac{t}{2k}$ of them are in $S(x,\ell+i)$. Indeed, let a vertex $y \in S(x,\ell)$
for $\ell \in \{ \lceil \frac{k}{2} \rceil,\lceil\frac{k}{2} \rceil+1,\cdots,k \}$. Then it is easy to see that the vertex $y$ has at most $cn^{k/2}$ neighbors in $S(x,j)$ for $j \in \{0,1,\cdots,\ell \}$, where $c$ is a constant that may depend on $k$. 

Therefore, there exists a set $S$ such that 
$$S \subset (A_1 \backslash A_0) \cap S(x,\ell+i) \cap \mathcal{N}(T \cap S(x,\ell))$$
and 
\begin{equation} \label{eq: |S|}
 \frac{(\frac{t}{2k})^2}{B(k)n^{(k-1)/2}}\leq |S| \leq \left(\frac{t}{2k} \right)^2,
 \end{equation}
where $B(k)$ is a large constant depending only on $k$. 

W.l.o.g, we can assume $x=0$. Let us count the maximum number of neighbors in $S(x,\ell)$ for a vertex $z$ in $S(0,\ell +i)$. Assume that $z=[1]^{\ell+i}[0]^{n-\ell-i}$. Let $y \in S(0,\ell)$ be  $y=[1]^a[0]^{\ell+i-a}[1]^{\ell-a}[0]^{n-2\ell-i+a}$. Then we have $\ell+i-a+\ell-a \leq k$ since $y$ is a neighbor of $z$. It is easy to see that $|\mathcal{N}(z) \cap S(0,\ell)|= \sum_{a=0}^{l+i} \binom{l+i}{a} \binom{n-\ell-i}{\ell-a}$. Therefore, as long as $n$ and $B(k)$ are large enough, (\ref{eq: |S|}) is valid. 

Now, let us consider the set $D: = \mathcal{N}(S) \cap S(x,\ell+i+k)$. Let $d$ be the number of edges between $S$ and $D$. We can express $D$ as a union of disjoint sets $D_i$, where $D = \cup_{i=1}^{\binom{\ell+i+k}{\ell+i}} D_i$. Each element in $D_i$ has $i$ neighbors in $S$. Let $|D_i|=d_i$ and $|S|=s$. Then we have 
$$d=\sum_{j=1}^{\binom{\ell+i+k}{\ell+i}}j d_j=s\binom{n-\ell-i}{k}.$$
Let $R_i=A_0 \cap D_i$ and $|R_i|=r_i$. We have
\begin{align}
        r:=\sum_{j=1}^{\binom{\ell+i+k}{\ell+i}}j r_j & \geq s\left(\frac{N}{2}-2t-\sum_{i=1}^{k-1}\binom{n}{i}\right) \\
        & = \left(\frac{n^k}{2k !}-\frac{1}{2}\left(\sum_{m=1}^{k-1}m \right)\frac{n^{k-1}}{k!}-\frac{1}{2}\frac{n^{k-1}}{(k-1)!} +O\left(\max\{n^{k-2},n^{k/2}\}\right)\right)s     \label{eq: r=jr_j}
\end{align}

and also
\begin{align}
\mathbb{E}(r) &=pd\\
&=p\binom{n-\ell-i}{k}s \\
&=\left(\frac{1}{2}-\left(\frac{\sqrt{kk!}}{2}+\epsilon \right) \frac{ \sqrt{\log n}}{n^{k/2}} \right)\binom{n-\ell-i}{k}s\\
&=\left(\frac{n^k}{2k!}- \frac{1}{2}\left(\sum_{m=0}^{k-1}(\ell+i+m)\right)\frac{n^{k-1}}{k!}+O\left( \max \{n^{k-2},n^{k/2} \sqrt{\log n}\}\right) \right)s. \label{eq: E(r)}
\end{align}

By comparing (\ref{eq: r=jr_j}) and (\ref{eq: E(r)}), we have that $r \geq \mathbb{E}(r)+\frac{n^{k-1}}{k!}s$ for $n$ large enough. Since the event $\{x \in A_3 \backslash A_2 \}$  implies the event $\{r \geq \mathbb{E}(r)+\frac{n^{k-1}}{k!}s\} ,$ we have 
\begin{align*}
    \mathbb{P}( \exists x \in A_3 \backslash A_2) & \leq 2^n \binom{n^k}{t} \binom{n^kt}{s}t^2\sum_{d_1,d_2,...d_{\binom{\ell+i+k}{k}}} \mathbb{P}\left(r \geq \mathbb{E}[r]+\frac{ n^{k-1}}{k!}s \right).
\end{align*}

Indeed, we have at most $2^n$ vertices in $Q_{k,n}$, $\binom{n^k}{t}$ choices for the set $T$, $\binom{n^kt}{s}$ choices for the set $S$ and at most $t^2$ values that $s$ can take.

Thus, from Lemma~\ref{lemmaA8}, there exists a constant $c >0$ such that 
$$\mathbb{P}\left(r \geq \mathbb{E}[r]+\frac{n^{k-1}}{k!}s \right) \leq \left(2\frac{n^{k-1}s}{k!} \right)^{\binom{\ell+i+k}{k}}\exp\left(- \frac{(2 \frac{n^{k-1}s}{k!})^2}{D\left(\binom{\ell+i+k}{k} \right)}\right).$$
Since $D\left(\binom{\ell+i+k}{k} \right)=\sum_{j=1}^{\binom{\ell+i+k}{k}} j^2d_j$ and $\sum_{j=1}^{\binom{\ell+i+k}{k}}j d_j=s\binom{n-\ell-i}{k}$, there exists a big absolute constant $c$ such that 
$$D\left(\binom{\ell+i+k}{k} \right)=s\binom{n-2k}{k}+\sum_{j=1}^{\binom{\ell+i+k}{k}}(j^2-j)d_j \leq cs\binom{n-2k}{k}.$$
\\
Thus, we have 
$$\mathbb{P}[\exists x \in A_3\backslash A_2] =o(1). \qedhere$$

\end{proof}

For the case where $k=2$, we will give a separate proof for $\mathbb{P}(A_3 = A_2) = 1 - o(1)$. 

Before stating the main lemma, we need to show two auxiliary lemmas. 

Let us first give some definitions.

Let $S \subset S([0]^n,4)=\binom{[n]}{4}$ and regard $S$ as 4-uniform hypergraph. Define $$d_S(i,j):=\{\{k,\ell\}:\{i,j,k,\ell\} \in S\}$$
and 
$$||S||_2:=\sum_{\{i,j\} \in \binom{[n]}{2}}\binom{d_S(i,j)}{2}.$$

Define 
$$d_S(i,j,k):=\{\ell:\{i,j,k,\ell\} \in S\}$$
and
$$||S||_3:=\sum_{\{i,j,k\} \in \binom{[n]}{3}}\binom{d_S(i,j,k)}{2}.$$

Consider $D=\mathcal{N}(S) \cap S(x,6)$. We can express $D$ as a union of disjoint sets $D_i$, where $D = \cup_{i=1}^{15} D_i$. Each element in $D_i$ has $i$ neighbors in $S$. Let $|D_i|=d_i$. This following lemma helps us establish a relation between $D(15)$ and $||S||_2$, $||S_3||$, which is very crucial in the proof of Lemma \ref{lemmaA18}. 

\begin{lemma} \label{lemma:S_2,S_3}
$$\sum_{i=2}^{15}\binom{i}{2}d_i=||S||_2 +(n-8)||S||_3$$
\end{lemma}
\begin{proof}
We will count the unordered triples $(A,B,W)$ in two ways, where $A,B \in S$, $A \neq B$, $W \in \binom{[n]}{6}$, and $A,B \subset W$. 

First fix a set $W \in D_i$, there are $\binom{i}{2}$ choices for $(A,B)$. Thus the total number of the unordered triples $(A,B,W)$ is $\sum_{i=2}^{15}\binom{i}{2}d_i$. 

Now fix $(A,B)$. Since $A,B$ are distinct $4$-sets and both of them are subset of the same $6$-set, either $|A \cap B|=2$ or $|A \cap B|=3$. 

If $|A \cap B|=2$, then $|A \cup B|=6$ and thus there is exactly one $6$-set containing both $A$ and $B$. 

If $|A \cap B|=3$, then $|A \cup B|=5$ and thus there are $n-5$ 6-sets containing both $A$ and $B$. 

Let $$P_2:=\{\{A,B\} \subset S:|A\cap B|=2\}$$ and 
$$P_3:=\{\{A,B\} \subset S: |A\cap B|=3\}.$$
Then, we have 
$$\sum_{i=2}^{15}\binom{i}{2}d_i=P_2 +(n-5)P_3.$$

Note that $||S||_2=P_2+3P_3$ and $||S||_3=P_3$. This completes the proof.
\end{proof}

We will need the following counting lemma, similar to the one in \cite{balogh_bollobas_morris_2009}, to estimate the number of choices of the set $S \in S(x,4)$ in the proof of Lemma ~\ref{lemmaA18}.
\begin{lemma} \label{lemma:size of clustered struture}
Let $G$ be a $4$-uniform hypergraph on $[n]$ and $H$ be a $4$-uniform sub-hypergraph of $G$. Assume that $e(G)=O(n^3)$  and $s=\Theta(n^2)$. We have 
$$ |\left\{H\subseteq G: e(H)=s,\; ||H||_2=m \right\} |
\leq \exp\left(\left( s-\frac{2m}{\binom{n-2}{2}} +o(s) \right)\log n \right).
$$
\end{lemma}
\begin{proof}
Call a pair \(P\in\binom{[n]}2\) large if $d_H(P)>\tau.$ We will specify $\tau$ later.  Let $K=\{P \in \binom{[n]}{2}: d_H(p) > \tau\}$ and write $ k:=|K|.$

Since every edge of \(H\) contains six pairs,
$$
\sum_{P\in\binom{[n]}2}d_H(P)=6s,
$$
and therefore
$$ k \leq \frac{6s}{\tau}.
$$

Let $L:=\{e \in E(H): e \supset P \; \text{for some } P \in K\}$ and $|L|=\ell$. For each $e \in L$, let $r_e$ be the number of large pairs contained in $e$. We have  
\begin{align}
\sum_{P\in K}d_H(P)&= \sum_{e \in L}r_e  \label{eq: double counting} \\
&=\ell+\sum_{e \in L}(r_e-1)\\
& \leq \ell+ \sum_{e \in L}\binom{r_e}{2}\\
&\leq \ell+(n-3)\binom{k}{2}, \label{eq:l+(n-3)}
\end{align}
where (\ref{eq: double counting}) is valid because both sides count $|\{(P,e): P \in K, e \in E(H), P \subset e \}|$ and (\ref{eq:l+(n-3)}) is valid because two distinct pairs are simultaneously contained in at most \(n-3\) different \(4\)-edges. 

Let $M=\binom{n-2}{2}.$ Since $d_H(P)\leq M$,
\begin{equation} \label{eq: p in K}
\sum_{P\in K}
\binom{d_H(P)}2
\leq
\frac{M}{2} \sum_{P\in K}d_H(P)\\
\leq
\frac{M}{2} \left(\ell+(n-3)\binom{k}{2}\right).
\end{equation}

For $P\notin K$, we have $d_H(P)\leq\tau$, and hence
\begin{equation}     \label{eq:P notin K}
\sum_{P\notin K} \binom{d_H(P)}2
\leq
\frac{\tau}{2} \sum_P d_H(P)\\
=3s\tau.
\end{equation}

Combining (\ref{eq: p in K}) and (\ref{eq:P notin K}), we have 
$$
m
\leq
\frac{M}{2}\left( \ell+(n-3)\binom{k}{2} \right)+3s\tau.
$$
Therefore
\begin{equation} \label{eq: l > 2m/M}
\ell \geq \frac{2\left(m-3s\tau-\frac{M}{2}(n-3)\binom{k}{2}\right)}{M}. 
\end{equation}

We now count the possible \(H\). To determine $H$ it is enough to choose $k$ large pairs in $K$, $\ell$ edges containing these pairs, and then the remaining $s-\ell$ edges. Thus,
$$|\left\{
H\subseteq G:
e(H)=s,\;
\|H\|_2=m
\right\} | \leq \sum_{k,\ell} \binom{n^2}{k} \binom{kn^2}{\ell} \binom{e(G)}{s-\ell},$$
where $k \leq \frac{6s}{\tau}$ and $ \frac{2\left(m-3s\tau-\frac{M}{2}(n-3)\binom{k}{2}\right)}{M}\leq \ell \leq s$.

Now $\tau:=\frac{M}{(\log n)^2}$ and thus $k \leq O((\log n)^2)$. We will use the fact that $(\frac{a}{x})^x$ with $a >0$ is increasing for $x \leq \frac{a}{e}$ and decreasing for $x \geq \frac{a}{e}$. Then we have for $c>0$
$$\binom{n^2}{k} \leq \left(\frac{en^2}{k}\right)^k \leq \left(\frac{en^2}{(c\log n)^2}\right)^{(c\log n)^2} \leq \exp (o(s\log n)).$$

Now let us bound the term $\binom{kn^2}{\ell}$. We will split it into two cases. If $\ell \geq \frac{s}{\log n}$, then for $c>0$
$$\frac{kn^2}{\ell} \leq \left(\frac{ekn^2}{\ell}\right)^{\ell} \leq\left( \frac{ekn^2}{cn^2} \right)^{cn^2}\leq \exp (n^2 O(\log \log n)) \leq \exp(o(s \log n)).$$
If $\ell \leq \frac{s}{\log n}$, we have 
$$\frac{kn^2}{\ell} \leq \left(\frac{ekn^2}{\ell}\right)^{\ell} \leq \exp (\ell O(\log n)) \leq \exp(O(s)) \leq \exp(o(s \log n)).$$

Now let us estimate $\binom{e(G)}{s-\ell}$. Since $e(G)=O(n^3)$, $s=\Theta(n^2),$ and $\ell \geq \frac{2m}{M}-o(s)$ we have 
$$
\binom{e(G)}{s-\ell} \leq \exp\left((s-\ell+o(s))\log n\right)
\leq \exp\left(\left(s-\frac{2m}{M}+o(s) \right) \log n\right). \qedhere
$$
\end{proof}

Now we are ready to state and prove the main lemma.
\begin{lemma} \label{lemmaA18}
 Consider the Boot($t$) process, where the initial infection probability $p= \frac{1}{2}-(1+\epsilon)\frac{ \sqrt{\log n}}{n}$ and $t=10n$. Then $A_3=A_2$ with high probability.
\end{lemma}
\begin{proof}
    Suppose $x \in A_3 \backslash A_2$.  Then we have 
$$|\mathcal{N}(x) \cap A_2| \geq \frac{N}{2}$$ and 
$$|\mathcal{N}(x) \cap A_1| \leq \frac{N}{2}-t.$$

Thus there exists a set $T \subset \mathcal{N}(x)$ with $|T|=t$ such that $T \subset A_2 \backslash A_1$ and 
$$|T \cap S(x,2)| \geq t-n,$$
since $|S(x,1)|=n$.

Let $y \in T$, and then we have 
$$|\mathcal{N}(y) \cap A_1| \geq \frac{N}{2}-t$$ and 
$$|\mathcal{N}(y) \cap A_0| \leq \frac{N}{2}-2t.$$
Thus every $y \in T$ must have at least $t$ neighbors in $A_1 \backslash A_0$.

For every $y \in (T \cap S(x,2))$ must have at least $t$ neighbors in $A_1 \backslash A_0$ among which at least $t-2(n-2)-(n-2)-2$ are in $S(x,4)$ since every vertex in $S(x,2)$ have $2(n-2)$ neighbors in $S(x,2)$, $n-2$ neighbors in $S(x,3)$ and $2$ neighbors in $S(x,1)$.Therefore, there exists a set $S$ such that 
$$S \subset (A_1 \backslash A_0) \cap S(x,4) \cap \mathcal{N}(T \cap S(x,2))$$
and 
$$ \frac{(t-n)(t-3(n-2)-2)}{6}\leq |S| \leq t^2$$
since every vertex in $S(x,4)$ has $6$ neighbors in $S(x,2)$.

Now consider $D=\mathcal{N}(S) \cap S(x,6)$. Let $d$ be the number of edges between $S$ and $D$. We can express $D$ as a union of disjoint sets $D_i$, where $D = \cup_{i=1}^{15} D_i$. Each element in $D_i$ has $i$ neighbors in $S$. Let $|D_i|=d_i$ and $|S|=s$. Then we have 
$$d=\sum_{i=1}^{15}i d_i=s\binom{n-4}{2}$$
Let $R_i=A_0 \cap D_i$ and $|R_i|=r_i$. We have

    $$r:=\sum_{i=1}^{15}i r_i \geq s\left(\frac{N}{2}-2t-4-\binom{4}{2}-(n-4)-\binom{4}{3}(n-4) \right)=s\left(\frac{N}{2}-2t-(5n-10) \right)$$
since every vertex in $S$ has at least $\frac{N}{2}-2t$ neighbors in $A_0$ and at most $4+\binom{4}{2}+(n-4)+\binom{4}{3}(n-4)$ neighbors in $A_0 \backslash S(x,6)$. 

We also have 
\begin{align*}
\mathbb{E}(r) &=pd\\
&=p\binom{n-4}{2}s \\
&=\left(\frac{1}{2}-\delta\right)\binom{n-4}{2}s
\end{align*}
where  $\delta=(1+\epsilon)\frac{ \sqrt{\log n}}{n}$.

Let $M=\binom{n-4}{2}$ and thus $N-M=5n-10$. Now we have $r-\mathbb{E}(r) \geq b_ns$, where $b_n=\delta M-2t-\frac{5n-10}{2}$. 

The main difficulty is that the number of possible choices for $S$ is too large if we use only the trivial upper bound $\binom{n^2t}{s}$. We therefore divide the analysis into two cases according to the size of $||S||_3$. When $||S||_3$ is large, we identify a subset $S'\subset S$ for which the number of possible choices is much smaller than $\binom{n^2t}{s}$. On the other hand, when $||S||_3$ is small, we apply Lemmas~\ref{lemma:S_2,S_3} and~\ref{lemma:size of clustered struture} to obtain an upper bound on the size of $S$.

\textbf{Case 1:} $||S||_3 \geq 3\eta ns$ with a sufficiently small $\eta >0.$ 

Let $\tau_3:=\eta n$ and define
$$
K_3
:=
\left\{
Q\in\binom{[n]}3:
d_S(Q)>\tau_3
\right\}.
$$

Since every edge of \(S\) contains exactly four triples, $\sum_{Q\in\binom{[n]}3}d_S(Q)=4s$, and therefore
$$ |K_3| \leq \frac{4s}{\tau_3}
= \frac{4s}{\eta n}
= O_\eta(n).
$$

The triples outside \(K_3\) contribute at most
$$
\begin{aligned}
\sum_{Q\notin K_3}
\binom{d_S(Q)}2
&\leq
\frac{\tau_3}{2}
\sum_{Q \notin K_3} d_S(Q)\\
&=2s\tau_3\\
&=2\eta ns.
\end{aligned}
$$

Thus, with $||S||_3 \geq 3\eta ns$, we have 
$$
\sum_{Q\in K_3} \binom{d_S(Q)}2 \geq \eta ns.
$$

Since $d_S(Q)\leq n-3$, we have $
\binom{d_S(Q)}2 \leq \frac{n}{2}d_S(Q).
$ Thus
$$
\sum_{Q\in K_3}d_S(Q)
\geq
2\eta s.
$$

Let $ L_3 :=
\left\{ e\in S: e\supseteq Q \text{ for some }Q\in K_3\right\}.$ Every edge of \(S\) contains at most four triples, and hence $ 4|L_3|
\geq \sum_{Q\in K_3}d_S(Q).$ Thus $$
|L_3|
\geq
\frac{\eta s}{2}.
$$

Therefore \(S\) contains a subset \(S'\) satisfying
$$
s':=|S'|
=\frac{\eta s}{2}
=
\Theta_\eta(n^2),
$$
such that every edge of \(S'\) contains a member of \(K_3\).

Now we need to count the number of choices for the set $S'$. There are at most
$$
\sum_{j\leq  cn}
\binom{\binom n3}{j}
= \exp(O_\eta(n\log n))
$$
possibilities for \(K_3\). Once \(K_3\) is fixed, the number of \(4\)-sets containing some member of \(K_3\) is at most
$
|K_3|(n-3)=O_\eta(n^2).
$

Hence the number of possibilities for \(S'\) when $K_3$ is fixed is at most
$$ \binom{cn^2}{s'} =\exp(O_\eta(n^2)) .
$$

The event $\{x \in A_3 \backslash A_2 \}$ implies the event $ \{r \geq \mathbb{E}(r)+\frac{\delta n^{2}}{4}s'\} $.
Therefore, we have 
$$ \mathbb{P}( \exists x \in A_3 \backslash A_2) \leq 2^n \binom{n^2}{t} \left(\sum_{j\leq  cn} \binom{\binom n3}{j} \right)\binom{cn}{s'}t^2\sum_{d_1,d_2,...d_{15}} \mathbb{P}\left(r \geq \mathbb{E}[r]+\frac{\delta n^2}{4}s' \right)$$
since we have at most $2^n$ vertices in the graph $Q_{k,n}$,$\binom{n^2}{t}$ choices for the set $T$, $\sum_{j\leq  cn}
\binom{\binom n3}{j}$ choices for the set $K_3$, $\binom{cn^2}{s'}$ choices for the set $S'$, and $t^2$ values that $s'$ can take.

Thus, from Lemma~\ref{lemmaA8}, we have 
$$\mathbb{P}\left (r \geq \mathbb{E}[r]+\frac{\delta n^2}{4}s' \right)\leq \left(\frac{\delta n^2}{2} \right)^{14} \exp \left(-\frac{(\delta n^2s')^2}{16D(15)} \right)$$

Since $D(15) =\sum_{i=1}^{15} i^2d_i$ and $\sum_{i=1}^{15}i d_i=s'\binom{n-4}{2}$, trivially  we have 
$D(15) \leq 15 s'n^2$.

Finally we have 
$$\mathbb{P}(\exists x \in A_3\backslash A_2)=o(1).$$ 
\medskip
\textbf{Case 2:} $||S||_3  <3\eta ns$

We now apply the Lemma ~\ref{lemma:size of clustered struture} to the hypergraph $G$ with vertex set $[n]$ and edge set $ E(G) = S(x,4)\cap \mathcal{N}(T\cap S(x,2))$. Note that $S$ is a sub-hypergraph of $G$.  Let $m_2:=||S||_2$ and $m_3:=||S||_3$. 

Every member of $T \cap S(x,2)$ is contained in at most $M_2=\binom{n-2}{2}$ members of \(S(x,4)\), and therefore
$$
e(G)\leq tM_2=O(n^3).
$$

Lemma ~\ref{lemma:size of clustered struture} implies that, once \(x,T,s,m_2\) are fixed, the number of possible choices of \(S\) is at most $$
\exp\left( \left(s-\frac{2m_2}{M_2}+o(1) \right)\log n \right).$$

By Lemma ~\ref{lemma:S_2,S_3}, we have 
\begin{align*}
D(15)
&:=\sum_{i=1}^{15}i^2d_i\\
&=
\sum_{i=1}^{15}i d_i
+
2\sum_{i=1}^{15}\binom{i}{2}d_i\\
&= Ms+2(m_2+(n-8)m_3).
\end{align*}

The event $\{x \in A_3 \backslash A_2 \}$ implies the event $ \{r \geq \mathbb{E}(r)+b_n s\} $. Then from Lemma ~\ref{lemmaA8} we have 
$$\mathbb{P}(\exists x \in A_3\setminus A_2)\leq 2^n \binom{n^2}{t}\exp\left( \left(s-\frac{2m_2}{M_2}+o(1) \right)\log n \right)t^2(2b_ns)^{14}
\exp\left(
-\frac{2b_n^2s^2}
{Ms+2(m_2+(n-8)m_3)}\right).
$$
Indeed, we have at most $2^n$ vertices in the graph $Q_{2,n}$,$\binom{n^2}{t}$ choices for the set $T$, $\exp\left( \left(s-\frac{2m_2}{M_2}+o(1) \right)\log n \right)$ choices for the set $S$, and $t^2$ values that $s$ can take.

Now we will try to bound this. 

Since $t=10n$ and $M=\binom{n-4}{2}$, for a sufficiently large $n$ we have $
\frac{b_n}{M}
\geq
\left(1+\frac{\epsilon}{2}\right)
\frac{\sqrt{\log n}}{n}.
$
In particular, if $
A:=\left(1+\frac{\varepsilon}{3}\right)^2,
$ then $
\frac{2b_n^2}{M}
\geq A\log n.
$

For all sufficiently large $n$, 
$$
\frac{2m_2}{Ms}+ \frac{2(n-8)m_3}{Ms}
\leq
u+12\eta+o(1),
$$
where $u=\frac{2m_2}{Ms}$. Note that $u \leq 1+o(1)$ due to (\ref{eq: l > 2m/M}). 

Therefore, we obtain
$$
\exp\left(-\frac{2b_n^2s^2}{Ms+2(m_2+(n-8)m_3)}\right)
=
\exp\left(-\frac{A+o(1)} {1+u+12\eta}s\log n \right).
$$

Now we have 
\begin{align*}
\exp\left(
-\frac{2b_n^2s^2}
{Ms+2(m_2+(n-8)m_3)}\right)\exp\left( \left(s-\frac{2m_2}{M_2}+o(1) \right)\log n \right) &=\\ 
\exp\left(\left[1-u- \frac{A}{1+u+12\eta} +o(1) \right] s\log n \right).
\end{align*}

We now choose $\eta>0$ sufficiently small depending only on $\epsilon$. Consider $ f(u):= 1-u- \frac{A}{1+u+12\eta}.$ Then $
f'(u) =-1+\frac{A}{(1+u+12\eta)^2}. $

If the maximum occurs in the interior, then $ 1+u+12\eta=\sqrt{A},$ and hence $$
f(u) = 2+12\eta-2\sqrt{A}.
$$

Thus, by taking \(\eta=\eta(\epsilon)>0\) sufficiently small,
$$
2+12\eta-2\sqrt A<0.
$$

At the endpoint \(u=0\),
$$
f(0)= 1-\frac{A}{1+12\eta}<0
$$
for a sufficiently small \(\eta\). 

For $u \geq 1$, it is obvious that $$f(u)<0.$$ Hence there exists a constant $c_\epsilon>0 $ such that uniformly over all \(u\),
$$
1-u-
\frac{A}{1+u+12\eta}
\leq
-c_\epsilon.
$$

Finally, it is easy to see that $2^n\binom{n^2}{t}t^2(2b_ns)^{14} \leq \exp(o(s\log n) )$

Therefore, $$\mathbb{P}(\exists x \in A_3 \setminus A_2) =o(1).\qedhere$$
\end{proof}

Finally, we are ready to prove the lower bound of Theorem~\ref{theorem:main}.

\begin{proof}
Consider the process Boot($t$), where the initial infection probability $p= \frac{1}{2}-\left(\frac{\sqrt{kk!}}{2}+\epsilon\right)\frac{\sqrt{\log n}}{n^{k/2}}$. If we can prove that $\mathbb{P}(A_2 = V(Q_{k,n})) = o(1)$ for Boot($t$), in conjunction with Lemmas ~\ref{lemmaA17} and \ref{lemmaA18}, it implies that the process Boot($t$) does not percolate with high probability. Consequently, this suggests that the process Boot($0$) does not percolate with high probability.

We will first show that for Boot($t$),
$$\mathbb{P}(x \in A_2)=\frac{1}{2}+o(1).$$

For a fixed vertex $x$, denote 
$$E_x:=|\mathcal{N}(x) \cap A_0|$$
and 
$$F_x:=|\mathcal{N}(x) \cap (A_1\setminus A_0)|.$$

Then $|\mathcal{N}(x) \cap A_1|=E_x+F_x$. 

We have $\mathbb{E}(E_x)=pN=\frac{1}{2}N-N\left(\frac{\sqrt{kk!}}{2}+\epsilon\right)\frac{\sqrt{\log n}}{n^{k/2}}.$ Let $\Delta_n'=N\left(\frac{\sqrt{kk!}}{2}+\epsilon\right)\frac{\sqrt{\log n}}{n^{k/2}}$. 

Consider the event 
$$E:=\left\{ E_x \leq \frac{N}{2}-\frac{\Delta_n'}{2}\right\}$$

By Lemma \ref{lemmaA1}, 
$$\mathbb{P}(E^c)= \mathbb{P}\left(E_x > \frac{N}{2}-\frac{\Delta_n'}{2}\right) \leq \exp\left(-c\frac{\Delta_n'^2}{N}\right)=o(1),$$ where $c$ is a constant.

Trivially, $\mathbb{P}(E^c | x \notin A_0)=\mathbb{P}(E^c)=o(1)$.

Now let us analyze the size of $F_x$. Let $y \in \mathcal{N}(x)$. We have 
\begin{align*}
\mathbb{P}\left(y \in (A_1 \setminus A_0)| x \notin A_0\right)&=\mathbb{P}(y \notin A_0,|\mathcal{N}(y) \cap A_0| \geq \frac{N}{2}-2t| x \notin A_0)\\
&=(1-p)\mathbb{P}(\text{Binomial}(N-1,p) \geq \frac{N}{2}-2t).
\end{align*}

By exactly the same analysis for $\tau_n$ in the proof of Lemma ~\ref{lemma:0.5+0.5delta}, we have 
$$\mathbb{P}(\text{Binomial}(N-1,p) \geq \frac{N}{2}-2t)=n^{-\alpha'+o(1)},$$ where $\alpha'=\frac{2(\frac{\sqrt{kk!}}{2}+\epsilon)^2}{k!}$. 

We also have 
$$\mathbb{E}(F_x| x \notin A_0) \leq N n^{-\alpha'+o(1)}=n^{k-\alpha'+o(1)}$$

By Markov's inequality, 
$$\mathbb{P}(F_x \geq \frac{\Delta_n'}{2}-2t | x \notin A_0) \leq \frac{\mathbb{E}(F_x | x \notin A_0)}{\Delta_n'}=o(1)$$
since $\alpha' > \frac{k}{2}$. 

Therefore, $$\mathbb{P}(x \notin A_2 | x \notin A_0)=1-o(1).$$ and thus $$\mathbb{P}(x \in A_2)= \frac{1}{2}+o(1).$$

From Lemmas \ref{lemma: graph partition} and \ref{lemma: graph partition for hypercube}, we can observe that there exists a partition of the set $V(Q_{k,n})$ as follows:
$$V(Q_{k,n})= \cup_{i=1}^m B_i,$$
where, for each $i \in [m]$, and for any $x$ and $y \in B_i$, we have $d_H(x,y) \geq 4k+1$. Furthermore, the number of subsets $m$ in this partition satisfies $m=\sum_{i=0}^{4k}\binom{n}{i}$.

Given the partition of $V(Q_{k,n})$, there exists a set $B_i$ such that $|B_i| \geq \frac{2^n}{n^{4k}}.$ Furthermore, for any $x$ and $y \in B_i$, it holds that $d_H(x,y) \geq 4k+1$. With the independence condition in place, we have 
\begin{align*}
\mathbb{P}(x \in A_2 \; \text{for all} \;  x \in B_i) &=\left(\frac{1}{2}+o(1)\right)^{|B_i|}\\
&= o(1). \qedhere
\end{align*}
\end{proof}

\section{Discussion and Open Problems}

Question 1. What is the third-order term of $p_c(Q_{k,n},\frac{N}{2})$ ?

In \cite{balogh_bollobas_morris_2009}, the third-order term of the critical probability was determined up to a constant factor. It is therefore natural to ask whether an analogous result can be obtained in our setting. At present, however, we do not even have a conjecture for the possible form of the third-order term.

Question 2. What is $p_c(Q_{k,n},r)$ for other ranges of $r$?

 When $r=2$ and $3$ this problem was solved by the same author in \cite{zhu2026bootstrappercolationgeneralizedhamming}. However, for other ranges of $r$ it remains widely open.

\section*{Acknowledgments}
The author would like to thank Prof.~Alexander Barg for suggesting this problem and providing valuable guidance. Additionally, the author extends thanks to
Prof.~J{\'o}zsef Balogh for answering questions related to \cite{balogh_bollobas_morris_2009}, Yihan Zhang for insightful discussions, and Anna Geisler for pointing out \cite{collares2024universalbehaviourmajoritybootstrap}.

\section*{AI Disclosure}
ChatGPT 5.6 Plus was used in the preparation of this paper. In particular, during a discussion concerning accurate approximations for binomial random variables in the moderate-deviation regime, ChatGPT suggested several relevant results, including Lemmas~\ref{lemma: Petrov} and~\ref{lemma:Gordon}. ChatGPT was also used to assist with calculations and to verify the computations appearing in Lemmas ~\ref{lemma: binom approx} and \ref{lemma:size of clustered struture}. All mathematical arguments and claims in this paper were developed by the author, who bears full responsibility for the correctness of the results.

\bibliographystyle{plain}
\bibliography{main.bib}
\end{document}